\documentclass [11pt]{article}
\usepackage{amsfonts,mathrsfs,graphicx,latexsym,bm}
\usepackage{amsmath,amsthm,latexsym,amssymb,amsxtra}
\usepackage{indentfirst}
\usepackage[driverfallback=dvipdfmx]{hyperref}
\usepackage{verbatim}

\makeatletter

\newtheorem{theorem}{Theorem}[section]{\bf}{\it}
\newtheorem{proposition}{Proposition}[section]{\bf}{\it}
\newtheorem{lemma}{Lemma}[section]{\bf}{\it}
\newtheorem{corollary}{Corollary}[section]{\bf}{\it}
\newtheorem{remark}{Remark}[section]{\it}{\rm}
\newtheorem{definition}{Definition}[section]{\it}{\rm}
{\it}{\rm}

\renewcommand{\sec}[1]{\setcounter{equation}{0}\section{#1}}
\newcommand{\bthm}{\begin{theorem}}
\newcommand{\ethm}{\end{theorem}}
\newcommand{\blem}{\begin{lemma}}
\newcommand{\elem}{\end{lemma}}
\newcommand{\bpp}{\begin{proposition}}
\newcommand{\epp}{\end{proposition}}
\newcommand{\bprf}{\begin{proof}}
\newcommand{\eprf}{\end{proof}}
\newcommand{\brem}{\begin{remark}}
\newcommand{\erem}{\end{remark}}
\newcommand{\bcor}{\begin{corollary}}
\newcommand{\ecor}{\end{corollary}}
\newcommand{\bdefi}{\begin{definition}}
\newcommand{\edefi}{\end{definition}}

\newcommand{\bprfp}[1]{\begin{proof}[Proof of Proposition \nopunct] \ref{#1}. \ \  }
\newcommand{\bprfl}[1]{\begin{proof}[Proof of Lemma \nopunct] \ref{#1}. \ \  }
\newcommand{\bprfc}[1]{\begin{proof}[Proof of Corollary \nopunct] \ref{#1}. \ \  }
\newcommand{\bprft}[1]{\begin{proof}[Proof of Theorem \nopunct] \ref{#1}. \ \  }

\newcommand{\lrc}{\lrcorner} 

\newcommand{\s}{\ }

\newcommand{\nr}{\\ && }

\newcommand{\nno}{ \nonumber }

\newcommand{\beqr}{\begin{eqnarray*} && }
\newcommand{\eeqr}{\end{eqnarray*}}
\newcommand{\beq}{\begin{eqnarray}}
\newcommand{\eeq}{\end{eqnarray}}
\newcommand{\beqs}{\begin{eqnarray*}}
\newcommand{\eeqs}{\end{eqnarray*}}
\newcommand{\bequ}{\begin{equation}}
\newcommand{\eequ}{\end{equation}}
\newcommand{\bequs}{\begin{equation*}}
\newcommand{\eequs}{\end{equation*}}
\newcommand{\bcase}{\begin{cases}}
\newcommand{\ecase}{\end{cases}}
\newcommand{\btab}{\begin{tabular}[t]}
\newcommand{\etab}{\end{tabular}}

\newcommand{\vphi}{\varphi}
\newcommand{\pa}{\partial}

\newcommand{\fr}{\frac}

\newcommand{\sr}{\stackrel}

\newcommand{\wg}{\wedge}

\newcommand{\na}{\nabla}

\newcommand{\dt}{\delta}
\newcommand{\ap}{\alpha}

\newcommand{\tta}{\theta}

\newcommand{\og}{\omega}
\newcommand{\Og}{\Omega}

\newcommand{\ot}{\otimes}

\newcommand{\subs}{\subset}

\newcommand{\bw}{\bigwedge}

\newcommand{\di}[1]{\langle #1 \rangle} 
\newcommand{\db}[1]{\{ #1 \}} 

\newcommand{\mbbr}{\mathbb{R}}

\begin{document}

\title{{$L^2$-estimates on $p$-convex Riemannian manifolds}
\footnotetext{Mathematics Classification Primary(2010):32F10, 53C21,53C23. \\ \hspace*{6mm}{Keywords: $p$-convex; Diederich-Forn{\ae}ss type exponent; uniform Carleman type estimate; finiteness theorem.}\\
\hspace*{6mm}{School of Mathematics and Key Laboratory of Mathematics for Nonlinear Sciences, Fudan
University, Shanghai 200433, China.}\\
\hspace*{5mm} {Email Address: qingchunji@fudan.edu.cn; xshliu@fudan.edu.cn; 12110180009@fudan.edu.cn.\\
\hspace*{6mm}{Partially
supported by NSFC11171069 and NSFC11322103.}}}}
\author{Qingchun Ji\ \ \ \ \ \ Xusheng Liu\ \ \ \ \ \ Guangsheng Yu }

\date{}

\maketitle

\begin{abstract}
 In this paper, we establish various $L^2$-estimates
 for the exterior differential operator on $p$-convex Riemannian manifolds
 in the sense of Harvey and Lawson. As applications, we establish a Carleman type estimate which is uniform with respect to both weight functions and domains, and we also obtain topological restrictions for
a Riemannian manifold to be $p$-convex.
\end{abstract}

\begin{center}
\tableofcontents
\end{center}

\setcounter{section}{-1}

\sec{Introduction}

The theory of convexity is a cornerstone of geometry, analysis and
related areas in mathematics. Recently in a series of
articles(\cite{HL1},\cite{HL2},\cite{HL3} and references therein),
Harvey and Lawson systematically explored the notions of
plurisubharmonicity and convexity in the context of differential
geometry. It has a long history for the concepts of pseudoconvexity
and plurisubharmonicity in several complex analysis and complex
geometry, but it has rare attention in more geometric situations
until Harvey and Lawson's innovative development in geometric
convexity. They also studied potential theory for geometric
plurisubharmonic functions and interesting applications to the
theory of nonlinear partial differential equations. A number of
results in complex analysis and complex geometry turn out to carry
over to more general setting. In \cite{HL3}, Harvey and Lawson
introduced the notion of $p$-convexity and $p$-plurisubharmonicity
on Riemannian manifolds. They obtained a deep result which is an
analogue of the Levi problem in complex analysis, i.e., local
$p$-convexity implies global $p$-convexity. This hopefully will
enrich the function theory in geometric analysis. For a compact
Rimannian manifold with smooth boundary, the concept of
$p$-convexity was first introduced by Sha(\cite{Sh}). In \cite{Sh},
it was proved that any Rimannian manifold with non-negative
sectional curvature and p-convex boundary has the homotopy type of a
CW-complex of dimension $<p$. This result was later strengthened by
Wu(\cite{W1}). Note that in \cite{Sh}, the $p$-convexity of a
Riemannian manifold $(M,ds^2)$ with boundary is equivalent to that
$\partial M$ is strictly $p$-convex in the sense of Harvey and
Lawson. The notion of $p$-convexity in the sense of Harvey and
Lawson is different from that introduced by Andreotti and
Grauert(cf.\cite{AG} and \cite{AV}) in the context of complex
analytic geometry which is defined by certain conditions on the number of negative or positive eigenvalues of the Levi form. 
The main difference is that the notion of $p$-convexity in the sense of Andreotti and
Grauert only depends on the underlying complex structure(which is used to define the Levi
form), while in Riemannian case the notion of $p$-convexity of Harvey and Lawson does depend
on the given metric, and this feature brings difficulties in
introducing complete metric because the $p$-convex property may not be preserved.\\

Since the $L^2$-method has many profound applications in several
complex analysis and complex
geometry(see\cite{Be4},\cite{D1},\cite{D2},\cite{H1},\cite{H2},\cite{K},\cite{OT},\cite{S1},\cite{S2}
and references therein), we will establish in the present paper
various $L^2$-estimates for the exterior differential operator on
$p$-convex Riemannian manifolds in the sense of Harvey and Lawson.
In many situations, the choice of weight functions and estimates for solutions in $L^2$-method are
crucial in applications(see, e.g., \cite{B4}, \cite{D1}, \cite{DF}, \cite{GHS}, \cite{K} and \cite{OT}). Hence we will make emphasis on several
different types of $L^2$-estimate. In \cite{H}, the author
considered the $\bar{\partial}$-problem on (weakly)$q$-pseudoconvex
domains in $\mathbb{C}^n$, but no effort was made to obtain good
estimates for solutions. The method developed here can be used to
establish estimates for $\bar{\partial}$-problem on
(weakly)$q$-pseudoconvex K\"{a}hler manifolds. To explain the
technique clearly, we will first prove $L^2$-estimates and existence
results in Euclidean spaces, and then we will show how the technique
still works on Riemannian manifolds. We also discuss geometric
applications of the $L^2$-method on $p$-convex Riemannian manifolds.
We only consider the problems of existence and interior
regularity(for the minimal $L^2$-solutions) in the present paper, we
plan to consider the problems of extension of closed forms, boundary regularity of minimal solutions and
more geometric applications in subsequent work.\\

This paper is arranged as follows. In section 1, we will recall
related notions of $p$-convexity and $p$-plurisubharmonicity in the
sense of Harvey and Lawson and prove some results of exterior
algebra which will be used later in our estimate. This section is
ended with a lemma concerning the choices of weight functions.
Section 2 is devoted to prove a theorem on the existence of certain
defining functions which shows that we also have the
Diederich-Forn{\ae}ss type exponent in this case. From this result,
we can reproduce a theorem due to Harvey and Lawson(\cite{HL3})
which says that boundary $p$-convexity implies $p$-convexity. In
section 3, we will establish the basic $L^2$-estimate and existence
theorem for the exterior differential operator on $p$-convex open
sets in $\mathbb{R}^n$. Based on the apriori estimate obtained in
section 3, we prove a Berndtsson type result in section 4.
This kind of estimate involves two $p$-plurisubharmonic weights with
opposite signs in the exponent. Such estimate for
$\bar{\partial}$-problem on pseudoconvex domains was originally
obtained by Donnelly and
Fefferman(see\cite{DF},\cite{Be1},\cite{B1},\cite{B2}). In section
5, we discuss the minimal $L^2$-solution and estimate for the
minimal $L^2$-solution with respect to a fixed weight function. In section 6,
we establishes an estimate by using non-plurisubharmonic weights,
the idea of our proof goes back to \cite{Be2}. In section 7, these
$L^2$-estimates obtained in sections 2-6 will be generalized to
$p$-convex Riemannian manifolds in the sense of Harvey and Lawson.
As geometric applications, we consider topological restrictions for
a Riemannian manifold to be $p$-convex in the last section. We will
prove vanishing and finiteness theorems for the de Rham cohomology
groups for $p$-convex Riemannian manifolds(without additional
curvature assumptions). A uniform estimate of Carleman type(Lemma
\ref{lem85}) plays an important role in establishing these results.
Following H\"{o}rmander's idea(\cite{H1}) and using a uniform G{\aa}rding type estimate, we prove this Carleman
type estimate which is uniform with respect to domains and weights. Lemma
\ref{lem85} is different from H\"{o}rmander's original estimate in
complex analytic case which was proved on a fixed domain. This
estimate allows us to prove, without using
the approximation theorem  for closed forms, a finiteness theorem for non-compact manifolds which are strictly $p$-convex at infinity(not only for relatively compact domains with  strictly $p$-convex boundary, and the underlying metric is not assumed to be complete). In fact, Lemma \ref{lem85}
applied to a fixed weight function and an increasing sequence of domains
gives the finiteness theorem(Theorem \ref{thm81}), by a similar
argument, Lemma \ref{lem85} applied to a fixed domain and an
increasing sequence of weight functions also gives the approximation theorem for
closed forms(Theorem \ref{thm82}).

\section{Preliminaries}

In this section, we will collect some facts on exterior algebra for
later use and recall the notions of $p$-convexity and
$p$-plurisubharmonicity in the sense of Harvey and Lawson
(\cite{HL1},\cite{HL2}, \cite{HL3}).

\par
\bigskip

Here and throughout this paper, the convention is adopted for summation over pairs of repeated indices.
Let $(V,\langle\cdot,\cdot\rangle)$ ba a $n$-dimensional Euclidean
space, we denote by $\db{e_1,\cdots,e_n}$ an orthonormal basis of
$(V,\langle\cdot,\cdot\rangle)$ and denote by
$\db{\omega^1,\cdots,\omega^n}$ its dual basis. For any multi-index $J=(j_1,\cdots,j_p)$, we set $\omega^J=\omega^{j_1}\wedge\cdots\wedge\omega^{j_p}$.

\bdefi A quadratic form $\tta=\tta_{ij}\og^i\ot \og^j\in
V^{*}\otimes V^{*}$ is called $p$-positive (semi-)definite if any
sum of $p$ eigenvalues of the symmetric matrix $(\tta_{ij})$ is
positive(nonnegative) where $1\leq p\leq n$. \edefi By using the
inner product $\langle\cdot,\cdot\rangle$, we identify the space of
symmetric endomorphisms of $V$ with the space of quadratic forms.
Then, a self-adjoint endomorphism $F$ is $p$-positive
definite(resp., semi-definite) if and only if for any $p$-plane
$W\subseteq V$, the $W$-trace ${\rm tr}_WF:={\rm tr}(F\mid_W)$ is
positive(resp., non-negative).

Denote by $\bw^p$ the linear space of $p$-forms on $V$. For any
quadratic form $\tta=\tta_{ij}\og^i\ot \og^j$, we introduce a
self-adjoint linear operator on $\bw^p$ by setting \beq
F_{\tta}= {\tta}_{jk}\og^k\wg e_j\lrc \label{eq12} \eeq where
$\lrcorner$ means the interior product. It follows directly from the
definition of $F_\theta$ that \beq\theta_{jk}
g_{jK}g_{kK}&=&(\theta_{jk}e_j\lrc g)_{K}\cdot (e_k\lrc g)_{K} \nno \\
&=& \langle \theta_{jk}e_j\lrc g,e_k \lrc g \rangle \nno\\
&=& \langle F_\theta g,g\rangle. \label{eq13} \eeq for any
$g=g_J\omega^J\in\bw^p$ where $K$ runs over all strictly increasing multi-indices of length $p-1$.

Now we compute the eigenvalues of $F_\theta$ in terms of those of $\theta$. Let us denote the eigenvalues of  $(\theta_{ij})$  by
$$\lambda_1\leq\cdots\leq\lambda_n,$$ after an orthogonal transformation,
we have $$F_{\theta}=\sum_{j=1}^n \lambda_j\og^j\wg e_j\lrc.$$
For any multi-index $J$ with $|J|=p,$ set \beq
\lambda_{J}=\sum_{j\in J}\lambda_j, \label{eq15} \eeq then we have
\begin{eqnarray*}
F_{\theta} \og^J &=& \sum_{j=1}^n \lambda_j \og^j \wg e_j\lrc \og^J\\
&=& \sum_{j=1}^n \lambda_j \og^j \wg\sum_{a=1}^p
(-1)^{a-1}\delta_{jj_a}\og^{j_1}\wg \cdots\wg
\widehat{\og^{j_a}}\wg\cdots\wg \og^{j_p}\\
&=& \sum_{a=1}^p\lambda_{j_a}\og^J =\lambda_{J}\og^J
\end{eqnarray*}
where $\delta_{jj_1}$ is the Kronecker delta and the circumflex over
a term means that it is to be omitted. Therefore, we have \beq {\rm
the \ set \ of \ eigenvalues \ of}\ F_{\theta} \ {\rm are \ given \
by} \ \{\lambda_{J}\ \mid \  |J|=p\}. \label{eq16} \eeq

Let $F: \bigwedge^p\rightarrow\bigwedge^p$ be a self-adjoint linear
map, we have the following orthogonal decomposition \beq
{\bigwedge}^p={\rm Ker}F \oplus {\rm Im}F, \label{eq17} \eeq which
implies that $F$ induces an isomorphism $F|_{{\rm Im}F}:
 {\rm Im}F\rightarrow{\rm Im}F$.
We can therefore define \beq F^{-1}:=(F|_{{\rm Im}F})^{-1}:{\rm
Im}F\rightarrow{\rm Im}F \label{eq18} \eeq for any self-adjoint
linear map $F$. Notice that $F$ itself is not required
to be invertible in the above definition.\\

The following lemma records the basic estimate concerning the self-adjoint operator $F_\theta$ defined by \eqref{eq12}. \blem \label{lem11} Let
$\theta=\theta_{ij}\og^i\otimes \og^j$ be a quadratic form. If
$\theta-\tau\otimes\tau$ is $p$-positive semi-definite where
$\tau=\tau_i\og^i$ is a $1$-form on $V$ and $1\leq p\leq n$, then
$$\tau\wedge\xi\in {\rm Im}F_\theta$$ for any $(p-1)$-form $\xi$ and
we have the following estimate
$$\langle F_\theta^{-1}f,\tau\wedge\xi\rangle\leq \langle
F_\theta^{-1}f,f\rangle^{\frac{1}{2}}|\xi|$$ for any $p$-form $f\in
{\rm Im}F_\theta$, in particular
$$\langle F_\theta^{-1}(\tau\wedge\xi),\tau\wedge\xi\rangle\leq
|\xi|^2.$$ \elem

\bprfl{lem11} By definition, we have $F_{\tau\otimes\tau}=\tau\wedge
X_\tau\lrcorner$ where $X_\tau:=\tau_ie_i$. Let $\eta,\tilde{\eta}$
be arbitrary $p$-forms, it is clear that
$$\langle F_{\tau\otimes\tau}\eta,\tilde{\eta}\rangle=\langle X_\tau\lrcorner\eta,X_\tau\lrcorner\tilde{\eta}\rangle.$$ Now we assume
$F_\theta\eta=0$, since the quadratic form $\theta-\tau\otimes\tau$
is $p$-positive semi-definite, we obtain
$$0=\langle
F_\theta\eta,\eta\rangle\geq|X_\tau\lrcorner\eta|^2$$which implies
$X_\tau\lrcorner\eta=0$. Therefore, we get
$$\langle\tau\wedge\xi,\eta\rangle=\langle\xi,X_\tau\lrcorner\eta\rangle=0.$$ Altogether, we have proved that
 $$\tau\wedge\xi\in({\rm Ker}F_\theta)^\perp={\rm Im}F_\theta.$$According to \eqref{eq18}, $F_\theta^{-1}(\tau\wedge\xi)$ is
 well-defined.

 Finally, we turn to the desired inequality. The Cauchy-Schwarz inequality
 gives \begin{eqnarray*}\langle F_\theta^{-1}(\tau\wedge\xi),\tau\wedge\xi\rangle&=&
 \langle X_\tau\lrcorner
 F_\theta^{-1}(\tau\wedge\xi),\xi\rangle\\&\leq&\langle X_\tau\lrcorner
 F_\theta^{-1}(\tau\wedge\xi),X_\tau\lrcorner
 F_\theta^{-1}(\tau\wedge\xi)\rangle^{\frac{1}{2}}|\xi|\\&=&\langle F_{\tau\otimes\tau}
 \circ F_\theta^{-1}(\tau\wedge\xi),F_\theta^{-1}(\tau\wedge\xi)\rangle^{\frac{1}{2}}|\xi|\\&
 \leq&\langle F_{\theta}
 \circ
 F_\theta^{-1}(\tau\wedge\xi),F_\theta^{-1}(\tau\wedge\xi)
 \rangle^{\frac{1}{2}}|\xi|\\&=&\langle\tau\wedge\xi,F_\theta^{-1}(\tau\wedge\xi)\rangle^{\frac{1}{2}}|\xi|.
  \end{eqnarray*}Dividing both sides by $\langle F_\theta^{-1}(\tau\wedge\xi),\tau\wedge\xi\rangle^{\frac{1}{2}}$ gives
  $$\langle F_\theta^{-1}(\tau\wedge\xi),\tau\wedge\xi\rangle^{\frac{1}{2}}\leq
|\xi|.$$This is the second inequality claimed in this lemma. Since $\langle F_\theta^{-1}\cdot,\cdot\rangle$ defines a
positive semi-definite bilinear form on ${\rm
Im}F_\theta\bigcap\bigwedge^p$, for any $p$-form $f\in {\rm
Im}F_\theta$, we have \beq \langle
F_\theta^{-1}f,\tau\wedge\xi\rangle &\leq& \langle
F_\theta^{-1}f,f\rangle^{\frac{1}{2}}\langle
F_\theta^{-1}(\tau\wedge\xi),\tau\wedge\xi\rangle^{\frac{1}{2}}\nno \\
&\leq& \langle
F_\theta^{-1}f,f\rangle^{\frac{1}{2}}|\xi|.\nno\eeq which implies the first inequality, and the proof is complete. \eprf

Now let us recall the notions of $p$-plurisubharmonicity and
$p$-convexity in the sense of Harvey and Lawson.

Let $(M,ds^2)$ be a $n$-dimensional oriented Riemannian manifold.
Let $\db{e_1,\cdots,e_n}$ be locally orthonormal frame with dual
coframe $\db{\omega^1,\cdots,\omega^n}$. With the Levi-Civita
connection $D$, the Hessian of a function $\varphi$ on $M$ is given
by $D^2\varphi(X,Y)=XY\varphi-D_XY\varphi$.

\bdefi A smooth function $\vphi$ defined on an open set $\Og\subs M$
is said to be $p$-plurisubharmonic if its Hessian $D^2\varphi$ is
$p$-positive semi-definite on $\Og$ and we call $\varphi$ strictly
$p$-plurisubharmonic if $D^2\varphi$ is $p$-positive definite on
$\Og$.\edefi

It is easy to see that for a K\"{a}hler manifold $(M,ds^2)$ the notion of $p$-plurisubharmonicity is defined by the Levi form of the given function which only depends on the underlying complex structure. In general case, it depends on the given  Riemannian metric.

In \cite{HL3}, it was proved that a smooth function $\vphi$ is
$p$-plurisubharmonic if and only if the restriction of $\varphi$ to
any $p$-dimensional minimal submanifold is subharmonic. In what follows, (strict) plurisubharmonicity means  (strict)$1$-plurisubharmonicity.\\

Given a smooth function $\vphi$, we denote
$$F_{\vphi}=F_{D^2\varphi}=\vphi_{jk} \og^k\wg e_j \lrc
$$
where $D^2\vphi:=\vphi_{ij}\og^i\ot\og^j$ is the Hessian of
$\varphi$.  It is easy to show that the operator $F_{\vphi}$, acting on differential forms, is exactly given by the difference of the Lie derivative and covariant derivative with respect to the gradient of $\varphi$(see Lemma \ref{lem71}). This observation will be useful to allow us to carry out Morrey's trick handling the boundary term. \\

Due to \eqref{eq16}, we have the following
criterion for $p$-plurisubharmonicity of a smooth function:
\emph{$\vphi$ is $p$-plurisubharmonic(resp., strictly
$p$-plurisubharmonic) on a domain $\Omega\subseteq M$ if and only if
$F_\vphi$(acting on $p$-forms) is positive semi-definite(resp.,
positive definite) at each point of $\Omega$.}\\

If $\vphi$ is \emph{strictly $p$-plurisubharmonic } on $\Omega$, by choosing $e_i$'s to be eigenvectors of $D^2\vphi(x)$ at a given point $x\in\Omega$, 
it follows from \eqref{eq16} that \beq
(F_\vphi^{-1}g)_{J}=\lambda_{J}^{-1}g_{J} \label{eq19} \eeq holds
for any $g=g_{J}\og^J \in\bigwedge^p$ and any given multi-indices
$J$ satisfying $|J|=p$ where $\lambda_{J}$ is defined by
\eqref{eq15} with $\theta=D^2\vphi(x)$. If the function $\vphi$ is further assumed to be
\emph{strictly plurisubharmonic}, we denote by
$(\vphi^{jk})$ the inverse
matrix of the  Hessian matrix
$(\vphi_{jk})$, then we have
\beq
\langle F_\vphi^{-1} g, g\rangle &=&
\lambda_{J}^{-1}|g_{J}|^2 \nno \\
&=&
(\sum_{j\in J}\lambda_j)^{-1}|g_J|^2 \nno \\
&\leq&
\frac{1}{p^2}\sum_{j\in
J}\lambda_j^{-1}|g_{J}|^2 \nno \\
&=& \frac{1}{p^2} \vphi^{jk} g_{jK} g_{kK} \label{eq110} \eeq
for arbitrary
$g=g_J\og^J\in\bigwedge^p$ where we have used \eqref{eq13} and \eqref{eq15} in the last equality.\\

\bdefi A Riemannian manifold $(M,ds^2)$ is called
(strictly)$p$-convex if it admits a smooth
(strictly)$p$-plurisubharmonic proper exhaustion function. It is
called strictly $p$-convex at infinity if it admits a proper
exhaustion function which is strictly $p$-plurisubharmonic outside a
compact subset of $M$. \edefi

Let $\Og\subs M$ be a compact domain with smooth boundary $\pa \Og$.
Let $II_{\pa \Og}(X,Y)=\di{D_XY,\nu}$ be the second fundamental form
of the boundary with respect to the inward pointing unit normal
vector $\nu$.

\bdefi The boundary $\pa\Og$ is said to be $p$-convex if ${\rm tr}_W
II_x\ge 0$ for any tangential $p$-plane $W\subseteq T_x(\pa\Og)$ and
any $x\in \pa\Og$. If the above inequality is strict for any
tangential $p$-plane $W$, $\pa\Og$ will said to be strict
$p$-convex.\edefi

The notion of boundary convexity can be described in terms of a
defining function as follows. Let $\rho$ be a defining function for
$\Og$, by \eqref{eq13} and \eqref{eq16}, we know that  $\pa\Og$ is
$p$-convex if and only if
$$\rho_{ij}g_{iK}g_{jK}\geq 0 $$
holds on $\pa \Og$ for every $p$-form $g =g_J\og^J$ which satisfies
$$\sum_{i=1}^{n}\rho_{i} g_{iK}=0$$ for all multi-indices $K$ with $|K|=p-1$. In \cite{HL3}, it was proved that \emph{ if the boundary $\pa\Og$ is
$p$-convex, then the domain $\Og$ is $p$-convex} (this also follows
from our Theorem \ref{thm21} below).\\

The following lemma is useful for choosing weight functions in
various applications of $L^2$-estimates.

\blem \label{lem81} Let $(M,ds^2)$ be a $n$-dimensional Riemannian
manifold and $\omega$ be a continuous function on $M$. We have the following conclusions: \\
(i)If $(M,ds^2)$ is strictly $p$-convex, then there is a strictly
$p$-plurisubharmonic proper exhaustion function $\varphi\in
C^\infty(M)$ such that $F_\varphi+\omega{\rm Id}$ is $p$-positive
definite on $M$.\\ (ii)If $(M,ds^2)$ is strictly $p$-convex at
infinity, then there is a $p$-plurisubharmonic proper exhaustion
function $\varphi\in C^\infty(M)$ such that $F_\varphi+\omega{\rm
Id}$ is $p$-positive definite outside some compact subset of $M$. In
particular, $(M,ds^2)$ is $p$-convex.\\ (iii)Let $\varphi\in
C^\infty(M)$ be a $p$-plurisubharmonic proper exhaustion function.
For any constant $c\in \mathbb{R}$ and $\eta\in L^2(M,{\rm Loc})$,
there is a function $\psi\in C^\infty(M)$ such that
$0\leq\psi-\varphi$ is $p$-plurisubharmonic, $\varphi\equiv\psi$
when $\varphi<c$ and $\int_M|\eta|^2e^{-\psi}<+\infty$.\elem

\bprfl{lem81}(i) Let us begin with any strictly $p$-plurisubharmonic
exhaustion function $\phi\in C^\infty(M)$. Set \beq
\Lambda_\phi:=\lambda_1+\cdots+\lambda_p\nno \eeq where
$\lambda_1\leq\cdots\leq\lambda_n$ are the eigenvalue functions of
the Hessian $D^2\phi$ with respect to the underlying metric $ds^2$, then we
know by definition that $\Lambda_\phi>0$.  Assume, without loss of
generality, $\inf_M\phi=0$ and denote\beq \Omega_\nu:=\{x\in M \
\mid \ \phi(x)<\nu\} \ {\rm for} \ \nu=1,2,\cdots.\nno\eeq Since the
functions $\Lambda_\phi>0$ and $\omega$ are both continuous on $M$,
one can always find for each $\nu=1,2,\cdots$ a positive constant
$\sigma_\nu$ such that \beq \label{eq81}\sigma_\nu
\Lambda_\phi+p\omega>0\ \ {\rm holds \ on} \
\Omega_{\nu+1}\setminus\Omega_{\nu}.\eeq Now we choose a function
$\kappa\in C^\infty[0,+\infty)$ such that
\beq\label{eq82}\kappa^{\prime}(t)>0,\kappa^{\prime\prime}(t)> 0 \ {\rm for} \
t\geq 0, \ \kappa^{\prime}(\nu)>\sigma_\nu \ {\rm for} \
\nu=1,2\cdots,\eeq and
\beq\label{eq83}\kappa^{\prime}(0)\inf_{\Omega_1}\Lambda_\phi+p\sup_{\Omega_1}\omega>0.\eeq
Set $\varphi=\kappa\circ\phi$, then
$D^2\varphi=\kappa^{\prime}\circ\phi\cdot
D^2\phi+\kappa^{\prime\prime}\circ\phi\cdot d\phi\otimes d\phi$ and
consequently we have
$$\Lambda_\varphi\geq\kappa^{\prime}\circ\phi\cdot\Lambda_\phi.$$ The
construction of $\kappa$ implies that $F_\varphi+\omega {\rm Id}$ is
$p$-positive definite on $M$. \\ (ii) In this case, the proof is a
slight modification of the proof given above and we will keep the
notations the same as above. By definition, we have a proper
exhaustion function $\phi$ and a compact subset $S\subseteq M$ such
that $\phi$ is strictly $p$-convex in $M\setminus S$. Without loss
of generality, we assume $S\subseteq\Omega_{\frac{1}{2}}$. Choose
$\chi\in C^\infty(\mathbb{R})$ such that
$$\chi^{\prime}(t)>0, \chi^{\prime\prime}(t)>0 \ {\rm for} \ t>\frac{1}{2} \ {\rm and} \
\chi(t)=0 \ {\rm for} \ t\leq \frac{1}{2}.$$ It is easy to see that
$\chi\circ\phi$ is a $p$-plurisubharmonic proper exhaustion function
and strictly $p$-plurisubharmonic outside
$\overline{\Omega_{\frac{1}{2}}}$, in particular, we have proved
that $(M,ds^2)$ is $p$-convex.

Let $\kappa\in C^\infty[0,\infty)$ be a function which satisfies
\eqref{eq81} and \eqref{eq82} with $\phi$ being replaced by
$\chi\circ\phi$, then it is easy to check that
$\varphi:=\kappa\circ\chi\circ\phi$ is a $p$-plurisubharmonic proper
exhaustion function and that $F_\varphi+\omega{\rm Id}$ is
$p$-positive definite outside $\Omega_1$(note that in this case, we
can not have \eqref{eq83} because $\inf_{\Omega_1}\Lambda_\phi$ is
not necessarily positive). \\ (iii) Choose a smooth function
$\gamma$ defined on $\mathbb{R}$ such that
$$\gamma^{\prime}(t)\geq 0,\gamma^{\prime\prime}(t)\geq 0 \ {\rm for} \ t\in \mathbb{R},  \gamma(t)\equiv0 \ {\rm for} \ t<c$$ and $$
\gamma(c+\nu)>\nu+\log\int_{\Omega_{c+\nu+1}}|\eta|^2 \ {\rm for} \
\nu=1,2,\cdots$$where $\Omega_{c+\nu+1}$'s are sub-level sets of
$\varphi$. Set $\phi=\gamma\circ\varphi$, then we know by definition
that $0\leq \phi$ is $p$-plurisubharmonic and \beq
\int_M|\eta|^2e^{-\phi}&=&\left(\int_{\Omega_{c+1}}+\sum_{\nu\geq
1}\int_{\Omega_{c+\nu+1}\setminus\Omega_{c+\nu}}\right)|\eta|^2e^{-\phi}\nno
\\ &\leq& \int_{\Omega_{c+1}}|\eta|^2e^{-\phi}+\sum_{\nu\geq
1}e^{-\gamma(c+\nu)}\int_{\Omega_{c+\nu+1}\setminus\Omega_{c+\nu}}|\eta|^2\nno
\\ &\leq& \int_{\Omega_{c+1}}|\eta|^2e^{-\phi}+\sum_{\nu\geq
1}e^{-\nu}<+\infty.\nno\eeq It is easy to see that
$\psi:=\varphi+\phi$ is a desired function. The proof is
complete.\eprf

\sec{The Diederich-Forn{\ae}ss type exponent}
In this section, we
prove a Diederich-Forn{\ae}ss type result on the defining function
for $p$-convex open set with smooth boundary.

 \bthm \label{thm21}
Let $\Omega\Subset\mbbr^n$ be a
 $p$-convex open set with smooth boundary and let $r\in
C^\infty(\overline{\Omega})$ be a defining function for $\Omega$.
Then for any strictly $p$-plurisubharmonic function $\vphi\in
C^\infty(\overline{\Omega})$, there exist constants $K>0, \eta_0\in
(0,1)$ such that
for any $\eta\in (0,\eta_0)$
the function $\rho:=-(-re^{-K\vphi})^\eta$
is strictly $p$-plurisubharmonic on $\Omega$.
\ethm

\bprft{thm21}It suffices to show that $\langle F_\rho g,g\rangle>0$
for any $0\neq g\in \bigwedge^p$.
\par
By direct computation, we obtain \beq \langle F_{\rho}g,g\rangle&=&
\eta (-r)^{\eta-2}e^{-K\eta \vphi} \Big[
Kr^2(\langle F_{\vphi}g,g\rangle-K\eta \langle\na\vphi\lrc g,\na\vphi\lrc g\rangle) \nonumber \\
&& +(-r)\langle F_rg,g\rangle
+(1-\eta)\langle\na r\lrc g,\na r\lrc g\rangle \nonumber \\
&& +2K\eta r \langle\na r\lrc g, \na\vphi\lrc g\rangle \Big]\label{eq211}
\eeq Throughout the proof, we denote by $A_1,A_2,\cdots$ various
constants which are independent of $\eta,K$.

Since the boundary of $\Omega$ is assumed to be smooth, for any sufficiently small $\varepsilon >0$ there is a smooth map
$\pi:N_\varepsilon\rightarrow\partial\Omega$ such that
\beq
{\rm dist}(x,\partial\Omega)=|x-\pi(x)|, \
\forall x\in N_\varepsilon
\label{eq212}
\eeq
where
$N_\varepsilon :=\{x\in\Omega \ | \
r(x)>-\varepsilon\}.$
As the function
$r\in C^\infty(\overline{\Omega})$
is a defining function for $\Omega$,
there exists a constant $A_1>0$
which only depends on $\varepsilon$
such that
\beq
{\rm dist}(x,\partial\Omega)\leq-A_1r(x), \ A_1\leq |\nabla r(x)|,
\  \forall x\in N_\varepsilon.
\label{eq213}
\eeq
For any $g\in \bigwedge^p,x\in
N_\varepsilon,$ we define $p$-forms $g_1(x), g_2(x)$ as follows:
$$g_1(x)=\frac{1}{|\nabla r(x)|^2}\nabla r(x)\lrcorner dr(x) \wedge
g$$ and $$ \ g_2(x)=\frac{1}{|\na r(x)|^2}dr(x)\wedge\nabla r(x)\lrcorner
g.$$ It is easy to see that $$g=g_1(x)+g_2(x), \
|g|^2=|g_1(x)|^2+|g_2(x)|^2 $$ and \beq \nabla r(x)\lrcorner
g_1(x)= 0,\ |g_2(x)|\leq\frac{1}{|\nabla r(x)|}|\nabla r(x)\lrcorner
g| \label{eq214} \eeq for every $x\in N_\varepsilon$. From
\eqref{eq212} and the first inequality in \eqref{eq213}, there is a
constant $A_2>0$ such that

\beq
|\langle F_r g_1,g_1\rangle(x)-\langle F_r
g_1,g_1\rangle(\pi(x))|&=&|\int_0^1\frac{d}{dt}
\langle F_r
g_1,g_1\rangle(tx+(1-t)\pi(x))dt|
\nno \\
&\leq&-A_2r(x)|g|^2 \ \ \ \ \ \ \ \ \ \ \ \ \ \ \ \ \ \ \ \ \ \ \ \
\ \ \ \ \ \label{eq215} \eeq holds for any $x\in N_\varepsilon$. By
using the identity in \eqref{eq214} and the definition of
$p$-convexity, we get
$$\langle F_rg_1,g_1\rangle(\pi(x))\geq 0, \
\forall x\in N_\varepsilon.$$ Therefore, for any $x\in N_\varepsilon$, the
following estimate follows from \eqref{eq215}
$$\langle F_rg_1,g_1\rangle(x)\geq A_2r(x)|g|^2.$$ Taking into account
of the inequality in
\eqref{eq214} and $|g_1(x)|\leq |g|$,
the above estimate implies that
\beq
\langle F_rg,g\rangle(x)\geq
A_2r(x)|g|^2-\frac{A_3}
{|\nabla  r(x)|}|\nabla r(x)\lrcorner
g|\cdot |g|
\label{eq216}
\eeq
holds for any $x\in N_\varepsilon$ where $A_3>0$
is another constant.

Since $\vphi$ is strictly $p$-plurisubharmonic on $\overline{\Omega},$
there is a constant $\sigma>0$ such that
\beq
\langle F_\vphi
g,g\rangle(x) - \eta K |\nabla \vphi(x)\lrcorner
g|^2\geq(\sigma-A_4\eta K)|g|^2
\label{eq217}
\eeq
holds for any $x\in
\Omega$ where $A_4:=\sup_{\Omega}|\nabla \vphi|^2$.
From \eqref{eq211}
and \eqref{eq217},
there exists a constant $A_5>0$ such that
\beq
\langle F_\rho g,g\rangle(x)\geq \eta(-r)^{\eta-2}e^{-\eta
K\vphi}\Big[Kr^2(x)(\sigma-\frac{\eta}{1-\eta}A_4
K)-A_5\Big]|g|^2
\label{eq218}
\eeq
holds for any $x\in\Omega$.

When $K>\frac{4A_5}{\sigma\varepsilon^2}$ and
$\eta\in(0,\frac{\sigma}{2A_4K+\sigma})$,
\eqref{eq218} implies that
\beq
\langle F_\rho g,g\rangle\geq \frac{1}{4}\eta(-r)^{\eta-2}e^{-\eta
K\vphi}K\varepsilon^2\sigma|g|^2
\label{eq219}
\eeq
holds on
$\Omega\setminus N_\varepsilon.$

Similarly, for any constants $\eta\in(0,\frac{\sigma}{2A_4K})$ and
$K>\frac{4}{\sigma}
(A_2+\frac{\sigma^2}{4A_4}+2A_6^2+\sigma^2),
A_6:=\frac{A_3}{2A_1},$
 from \eqref{eq211}, \eqref{eq216} and \eqref{eq217}
it follows that
the following inequality holds on $N_\varepsilon$
\beq
\langle F_\rho g,g\rangle &\geq & \eta(-r)^{\eta-2}e^{-\eta K\vphi}\Big\{[K(\sigma-A_4\eta K)-A_2]r^2|g|^2 \nno \\
& &+ 2(A_6+A_4\eta K)|\nabla r\lrcorner
g|r|g| \nno \\
& &+(1-\eta)|\nabla r\lrcorner g|^2\Big\}\nno \\
&\geq & \eta(-r)^{\eta}e^{-\eta K\vphi}\Big[K(\sigma-A_4\eta
K)-A_2-\frac{2A_6^2+2A_4^2\eta^2K^2}
{1-\eta}\Big]|g|^2\nno \\
&\geq &\eta(-r)^{\eta}e^{-\eta K\vphi}
\Big(\frac{K\sigma}{2}-A_2-4A_6^2-\sigma^2\Big)|g|^2 \nno \\ &\geq
&\frac{1}{4}\eta(-r)^{\eta} e^{-\eta K\vphi}K\sigma|g|^2\nno \\
&=&\frac{K\eta\sigma}{4}(-\rho)|g|^2\label{eq220} \eeq By combining
\eqref{eq219} and \eqref{eq220}, we know Theorem \ref{thm21} is true
for any constant
$K>\frac{4}{\sigma}(A_2+\frac{\sigma^2}{4A_4}+\frac{A_5}{\varepsilon^2}+2A_6^2+\sigma^2)$
and $\eta_0:=\frac{\sigma}{2A_4K+\sigma}$.\eprf

\brem\label{rem21} (i) The constant $\eta$ is an analogue of the
Diederich- Forn{\ae}ss exponent in several complex variables(see
\cite{DiFo}). \\  (ii) By Theorem 3.1, we know that $\psi:=-\log(-\rho)$ is a
strictly $p$-plurisubharmonic proper exhaustion function on
$\Omega$, and this implies Theorem 3.10 in \cite{HL3}. \erem

\sec{The $L^2$-existence theorem}

Let $\Omega\subseteq \mathbb{R}^n$ be an open subset, $\varphi\in
C^1(\Omega)$. Following \cite{H1}, the weighted $L^2$-Hermitian
inner product of $p$-forms will be denoted by
$(\cdot,\cdot)_\varphi$ and the corresponding Hilbert space will be
denoted by $L_p^2(\Omega,\varphi).$ We will still denote by $d$ the
maximal(weak) differential operator(from $L_{p-1}^2(\Omega,\varphi)$
to $L_p^2(\Omega,\varphi)$ ) of the exterior differential. It is
easy to see that the formal adjoint of $d$ with respect to the weight
$\varphi$ is given by $\delta_\varphi:=e^{\vphi}\circ\dt\circ
e^{-\vphi}$ where $\delta$ is the codifferential operator on
$\mathbb{R}^n$. If $\Omega\Subset\mathbb{R}^n$ has smooth boundary
and $\varphi\in C^1(\overline{\Omega})$, integration by parts shows
that $C^\infty_p(\overline{\Omega})\cap{\rm Dom}(d^*_\varphi)=\{g\in
C^\infty_p(\overline{\Omega}) \ | \ \nabla\rho\lrcorner g=0 \ {\rm
on} \ \partial{\Omega}\}$ where $d^*_\varphi$ is the Hilbert space
adjoint of $d$ with respect to the weight $\varphi$ and $\rho$ is a defining
function of $\Omega$.

The following Kohn-Morrey-H\"{o}rmander type identity is crucial in
establishing basic apriori estimate.

\bpp\label{pp31} Let $\Og\Subset \mbbr^n$ be a domain with smooth
boundary. Assume that the defining function satisfies $|\na\rho|=1$
when restricted  to $\pa\Og$. Then we have the following identity:
\beq \|dg\|^2_{\vphi}+\|\dt_{\vphi}g\|^2_{\vphi} &=&\int_{\Og}
|\pa_jg_I|^2 e^{-\vphi}
+\int_{\Og} \langle F_{\vphi}g,g\rangle  e^{-\vphi} \nno \\
&& +\int_{\pa\Og} \langle F_{\rho}g,g\rangle e^{-\vphi} \label{eq31}
\eeq for $g\in C_p^\infty(\overline{\Omega})\cap {\rm
Dom}(d_\varphi^{*})(1\leq p\leq n)$. \epp 
\bprfp{pp31} Let
$\dt_k=e^{\vphi}\pa_k(e^{-\vphi}\cdot)$, then it is easy to see that
$$[\dt_k, \pa_j]=\pa_{j}\pa_{k}\vphi$$  holds on functions. By
definition, we have the following equalities \beqr
|dg|^2=|\pa_jg_I|^2 - \pa_jg_{kK}\pa_kg_{jK},\nr |\dt_{\vphi} g|^2=
\dt_jg_{jK}\dt_kg_{kK}. \eeqr Repeated use of the formula
$$\int_\Og \pa_jv w e^{-\vphi}
=\int_\Og -v\dt_j w e^{-\vphi} +\int_{\pa\Og}\pa_j\rho vw
e^{-\vphi}$$gives that \beq
\int_{\Og}\pa_jg_{kK}\pa_kg_{jK}e^{-\vphi}
&=&-\int_{\Og}g_{kK}\dt_{j}\pa_kg_{jK}e^{-\vphi}
+\int_{\pa\Og}g_{kK}\pa_kg_{jK}\pa_j\rho e^{-\vphi}\nno\\
&=&-\int_{\Og}g_{kK}\left(\pa_k\dt_jg_{jK}+[\dt_j,\pa_k]g_{jK}\right)
e^{-\vphi}
+\int_{\pa\Og}g_{kK}\pa_kg_{jK}\pa_j\rho e^{-\vphi}\nno \\
&=&\int_{\Og}\dt_kg_{kK}\dt_jg_{jK}e^{-\vphi}
-\int_{\pa\Og}g_{kK}\dt_jg_{jK}\pa_k\rho e^{-\vphi}\nno \\ &&
-\int_{\Og} g_{kK} \pa_{j}\pa_{k}\vphi g_{jK} e^{-\vphi}
+\int_{\pa\Og}g_{kK}\pa_kg_{jK}\pa_j\rho e^{-\vphi}.\nno \eeq From the
boundary  condition $$\pa_k\rho g_{kK}=0 \  {\rm on} \ \partial\Omega,$$ we know that, for any
fixed $K$ with $|K|=p-1$, $$g_{jK}\fr{\pa}{\pa x^j}{\rm \ defines \ a
\ tangent\ vector \ field \ of } \ \partial\Omega.$$ Consequently, we obtain \beq
0=\sum_{k=1}^n g_{kK}\pa_k ( \sum_{j=1}^n g_{jK} \pa_j\rho)= \sum_{j,k=1}^n
(g_{kK}\pa_kg_{jK}\pa_j\rho+\pa_{j}\pa_{k}\rho g_{kK} g_{jK}) \s
{\rm on} \s \pa\Og.\nno \eeq Therefore, we get \beq
 \|dg\|^2_{\vphi}+\|\dt_{\vphi}g\|^2_{\vphi}\nno
&=&\int_{\Og}|\pa_jg_I|^2 e^{-\vphi} +\int_{\Og}
\pa_{j}\pa_{k}\vphi g_{jK}g_{kK}  e^{-\vphi}\nno \\
&&+\int_{\pa\Og}\pa_k\rho g_{kK} \dt_jg_{jK} e^{-\vphi}
-\int_{\pa\Og}g_{kK}\pa_kg_{jK}\pa_j\rho e^{-\vphi}\nno \\
&=&\int_{\Og}|\pa_jg_I|^2 e^{-\vphi} +\int_{\Og} \pa_{j}\pa_{k}\vphi
g_{jK}g_{kK}  e^{-\vphi}
+\int_{\pa\Og}\pa_{j}\pa_{k}\rho  g_{jK}g_{kK} e^{-\vphi}\nno \\
&=&\int_{\Og}|\pa_jg_I|^2 e^{-\vphi} +\int_{\Og} \langle
F_{\vphi}g,g\rangle e^{-\vphi} +\int_{\pa\Og} \langle
F_{\rho}g,g\rangle e^{-\vphi} \nno\eeq which gives the desired
identity \eqref{eq31}. \eprf

To establish $L^2$-existence theorem, we also need the following
basic lemma from functional analysis due to H\"{o}rmander:
\blem\label{lem31}(Theorem 1.1.4 in \cite{H1}) Let
$H_1\sr{T}{\to}H_2\sr{S}{\to}H_3$ be a complex of closed and densely
defined operators between Hilbert spaces and let $L\subseteq H_2$ be
a closed subspace which contains ${\rm Im}(T)$. For any $f\in L\cap
{\rm Ker}(S)$ and any constant $C>0$,
the following conditions are equivalent\\
1.there exists some $u\in H_1$ such that
$Tu=f$ and $\|u\|_{H_1}\le C$.\\
2.$|(f,g)_{H_2}|^2\le C^2(\|T^{*}g\|_{H_1}^2+\|Sg\|_{H_3}^2)$ holds
for any $g\in L\cap{\rm Dom}(T^{*})\cap {\rm Dom}(S)$. \elem Now we
can prove a $L^2$-existence result for the exterior differential
operator.

\bthm\label{thm31} Let $\Og\subset \mbbr^n$ be a p-convex domain and
$\vphi\in C^2(\Omega)$ be a p-plurisubharmonic function over $\Og$.
Then for any closed p-form $f\in L^2_p(\Og,{\rm Loc})$ satisfying
$$\int_{\Og}\langle F_{\vphi}^{-1}f,f \rangle e^{-\vphi}<\infty$$
there exists a $(p-1)$-form
 $u\in L_{p-1}^2(\Og,\vphi)$ such that
$$
du=f, \|u\|^2_{\vphi}\le \int_{\Og} \langle F_{\vphi}^{-1}f,f
\rangle e^{-\vphi}
$$where $F_\varphi^{-1}$ is defined by \eqref{eq18} and it is assumed
implicitly that $F_{\vphi}^{-1}f$ is defined almost everywhere in
$\Omega$. \ethm 
\bprft{thm31} First we suppose
$\Og\Subset\mathbb{R}^n$ has smooth $p$-convex boundary. Then we
have, in formula \eqref{eq31}, $\langle F_{\rho}g,g \rangle \ge 0$
on $\pa \Og$, which implies \beq\label{eq32}
\|dg\|^2_{\vphi}+\|\dt_{\vphi}g\|^2_{\vphi} \ge \int_{\Og} \langle
F_{\vphi}g,g \rangle  e^{-\vphi} \eeq holds for any $g\in
C_p^\infty(\overline{\Omega})\cap {\rm Dom}(d_\varphi^{*})$. By
H\"{o}mander's density Lemma(\cite{H1}\cite{H2}), the above estimate
\eqref{eq32} holds for any $g\in{\rm Dom}(d_\varphi^{*})\cap{\rm
Dom}(d)$.

We will apply the Lemma \ref{lem31} to\\
$$
H_1=L^2_{p-1}(\Og,\vphi),
H_2=L^2_p(\Og,\vphi),
H_3=L^2_{p+1}(\Og,\vphi)
$$
and $S,T$ both given by the maximal differential operators of
exterior differentials. Since the $\langle F_{\vphi}\cdot,\cdot
\rangle $ is positive semi-definite, it follows from schwarz
inequality that \beq |\int_{\Og}\langle f,g\rangle e^{-\vphi}|^2
&=&|\int_{\Og} \langle F_{\vphi}F^{-1}_{\vphi}f,g \rangle
e^{-\vphi}|^2\nno \\ &\le&
(F_{\vphi}F^{-1}_{\vphi}f,F^{-1}_{\vphi}ff)_{\vphi}
(F_{\vphi}g,g)_{\vphi}\nno \\ &=& (F^{-1}_{\vphi}f,f)_{\vphi}
(F_{\vphi}g,g)_{\vphi}\nno \\  &\le& (F^{-1}_{\vphi}f,f)_{\vphi}
(\|T^{*}g\|^2_{H_1}+\|Sg\|^2_{H_2}). \nno \eeq Now from the Lemma
\ref{lem31}, it follows that there is a $(p-1)$-form $u\in
L^2_{p-1}(\Og,\vphi)$ such that
$$
du=f, \|u\|^2_{\vphi}\le \int_{\Og}\langle F_{\vphi}^{-1}f,f\rangle
e^{-\vphi}.
$$

For general case, by Theorem 3.4 in \cite{HL3}, there exists a
sequence of domains $\Omega_\nu(\nu=1,2,\cdots)$ with smooth
$p-$convex boundary such that $\Omega=\displaystyle{\cup_{\nu\geq
1}}\Omega_\nu$. For each $\nu\geq 1$, we have a solution $u_\nu\in
L_{p-1}^2(\Omega_\nu,\varphi)$ of $du_\nu=f$ such that
$$\int_{\Og_\nu} |u_\nu|^2e^{-\vphi}\le \int_{\Og_\nu} \langle F_{\vphi}^{-1}f,f \rangle
e^{-\vphi}.$$ By the estimate on $u_\nu$ we obtain the desired
solution by taking weak limit. The proof is complete.\eprf

Starting from any strictly $p$-plurisubharmonic proper exhaustion
function $\varphi$ and then using Lemma \ref{lem81} (iii), we have
the following corollary of Theorem \ref{thm31}.

\bcor \label{cor31} Let $\Og\subset \mbbr^n$ be a p-convex domain.
For any closed p-form $f\in L^2_p(\Og,{\rm Loc})$, there exists a
$(p-1)$-form
 $u\in L_{p-1}^2(\Og,{\rm Loc})$ such that
$ du=f$. \ecor

\sec{A Berndtsson type estimate} In this section, we will
establish a Berndtsson type result which involves two
$p$-plurisubharmonic weights with opposite signs in the exponent.
This kind of estimates for $\bar{\partial}$-problem was first
obtained by Berndtsson(see
\cite{Be1},\cite{Be3},\cite{B1}, \cite{B2} and references therein) and had its root in  Donnelly-Fefferman estimate(\cite{DF}).
The key for our proof is to establish the following apriori estimate.

$$\|\delta_{\varphi+\sigma\psi}g\|^2_{\varphi+\psi}
+\|dg\|^2_{\varphi+\psi}\geq \sigma^2\int_\Omega\langle F_{\psi}
g,g\rangle e^{-\varphi -\psi},\eqno{(*)}$$ for any $g\in {\rm
Dom}(d^{*})\cap C^\infty_p(\overline{\Omega})$ where $\varphi\in
C^\infty(\overline{\Omega})$ is a $p$-plurisubharmonic function,
$\psi\in C^\infty(\overline{\Omega})$ with $-e^{-\psi}$ being
$p$-plurisubharmonic and $\sigma\in (0,\frac{1}{2}]$ is a constant. The following proof involves a useful proceedure to introduce a twist factor into a 
known apriori estimate(see also \cite{Be3}, \cite{J}, \cite{S2}).

\bthm\label{thm41} \emph{Let $\Omega$ be a $p$-convex domain in
$\mbbr^n(1\leq p\leq n)$ and let $\varphi$ be a $p$-plurisubharmonic
function on $\Omega$, $\psi\in C^2(\Omega)$ be a function such that
$-e^{-\psi}$ is $p$-plurisubharmonic. For any constant $\ap\in
[0,1)$ and any closed $p$-form $f\in L_{p}^2(\Omega,{\rm Loc})$, if
$$\int_\Omega\langle F_\psi^{-1}f,f\rangle
e^{-\varphi+\ap\psi}<\infty$$
then there exists a $(p-1)$-form
$u\in L_{p-1}^2(\Omega,\varphi-\ap\psi)$ such that
$$du=f,\ \
\|u\|_{\varphi-\ap\psi}^2\leq\frac{4}{(1-\ap)^2}\int_\Omega\langle
F_\psi^{-1}f,f\rangle e^{-\varphi+\ap\psi},$$where $F_\psi^{-1}$ is
defined by \eqref{eq18} and it is required implicitly that
$F_\psi^{-1}f$ is defined almost everywhere in $\Omega$.} \ethm
\bprft{thm41} \ \  We consider first the case where $\Omega$ is a
bounded domain in $\mbbr^n$ with smooth boundary and
$\varphi,\psi\in C^\infty(\overline{\Omega})$.  \par We will apply
the Lemma \ref{lem31} to following weighted $L^2$-spaces of
differential forms
$$H_1=L^2_{p-1}(\Omega,\varphi+\frac{1-\ap}{2}\psi),
H_2=L^2_{p}(\Omega,\varphi+\frac{1-\ap}{2}\psi),
H_3=L^2_{p+1}(\Omega,\varphi+\frac{1-\ap}{2}\psi)$$ and
$$T=d\circ
e^{-\frac{1+\ap}{4}\psi},S=e^{-\frac{1+\ap}{4}\psi}\circ d.$$ In
order to use the Lemma \ref{lem31}, we need to show that the
following estimate \beq |(f,g)_{\varphi+\frac{1-\ap}{2}\psi}|^2
\leq\frac{4(F_\psi^{-1}f,f)_{\varphi-\ap\psi}}
{(1-\ap)^2}(\|e^{-\frac{1+\ap}{4}\psi}
\dt_{\varphi+\frac{1-\ap}{2}\psi}g\|^2_{\varphi+\frac{1-\ap}{2}\psi}
+\|e^{-\frac{1+\ap}{4}\psi}dg\|^2_{\varphi+\frac{1-\ap}{2}\psi})\eeq
holds for arbitrary
$g\in {\rm Dom}(d^{*})\cap C^\infty_{p}(\overline{\Omega}).$\\

Let $g\in {\rm Dom}(d^{*})\cap C^\infty_{p}(\overline{\Omega})$,
then the basic estimate with $\phi=\varphi+\psi$ gives \beq
\|dg\|_{\varphi+\psi}^2+\|\dt_{\varphi+\psi} g\|_{\varphi+\psi}^2
 \geq \int_\Omega\langle F_{\varphi+\psi} g,g\rangle e^{-\varphi-\psi}.\label{eq41}
\eeq
Since $$\dt_{\varphi+\psi}
g=\dt_{\varphi+\frac{1-\ap}{2}\psi}
g+\frac{1+\ap}{2}\nabla\psi\lrcorner g,$$
it follows that
$$\|\dt_{\varphi+\psi}
g\|_{\varphi+\psi}^2\leq
\frac{1+\epsilon}{\epsilon}\|\dt_{\varphi+\frac{1-\ap}{2}\psi}
g\|_{\varphi+\psi}^2+\frac{(1+\epsilon)(1+\ap)^2}{4}\|\nabla\psi\lrcorner
g\|_{\varphi+\psi}^2$$
for any positive constant $\epsilon.$\\
By
choosing $$\epsilon=\frac{1-\ap}{1+\ap},$$
the above inequality  becomes
\beq
\|\dt_{\varphi+\psi}
g\|_{\varphi+\psi}^2\leq
\frac{2}{1-\ap}\|\dt_{\varphi+\frac{1-\ap}{2}\psi}
g\|_{\varphi+\psi}^2+\frac{1+\ap}{2}\|\nabla\psi\lrcorner
g\|_{\varphi+\psi}^2.\label{eq42}
\eeq

Since $-e^{-\psi}$ is $p$-plurisubharmonic and
\begin{eqnarray*}
F_{-e^{-\psi}}&=&e^{-\psi}(\psi_{jk}
-\psi_j\psi_{k})dx^k\wedge\frac{\partial}{\partial
x^j}\lrcorner\\
&=&
e^{-\psi}(F_\psi-d\psi\wedge\nabla \psi\lrcorner),
\end{eqnarray*}
we know that $F_\psi-d\psi\wedge\nabla\psi\lrcorner$
defines a positive semi-definite operator on the space of
$p$-forms. This implies
\beq
\int_\Omega\langle F_{\psi} g,g\rangle e^{-\varphi-\psi}\geq\|
\nabla\psi\lrcorner g\|_{\varphi+\psi}^2.
\label{eq43}
\eeq
Substituting \eqref{eq42}, \eqref{eq43} into
\eqref{eq41},
the $p$-plurisubharmonicity of $\varphi$ gives
$$\frac{2}{1-\ap}\|\dt_{\varphi+\frac{1-\ap}{2}\psi}
g\|_{\varphi+\psi}^2 + \|d g\|_{\varphi+\psi}^2
 \geq \frac{1-\ap}{2}\int_\Omega\langle F_{\psi} g,g\rangle e^{-\varphi-\psi}
$$which further implies the desired Donnelly-Fefferman
type estimate $(*)$ with the constant $\sigma=\frac{1-\ap}{2}$ as
follows
\begin{eqnarray*}
\|e^{-\frac{1+\ap}{4}\psi}\dt_{\varphi+\frac{1-\ap}{2}\psi}g\|^2_{\varphi+\frac{1-\ap}{2}\psi}
+\|e^{-\frac{1+\ap}{4}\psi}dg\|^2_{\varphi+\frac{1-\ap}{2}\psi}
&=&
\|\dt_{\varphi+\frac{1-\ap}{2}\psi}g\|^2_{\varphi+\psi}
+\|dg\|^2_{\varphi+\psi}\\ &\geq&
\|\dt_{\varphi+\frac{1-\ap}{2}\psi}g\|^2_{\varphi+\psi}
+\frac{1-\ap}{2}\|dg\|^2_{\varphi+\psi}\\ &\geq&
\frac{(1-\ap)^2}{4}\int_\Omega\langle F_{\psi} g,g\rangle
e^{-\varphi -\psi}.
\end{eqnarray*}
Since $\psi$ is $p$-plurisubharmonic, the Cauchy-Schwarz inequality applied to the positive semi-definite Hermitian form
$(F_\psi\cdot ,
\cdot)_{\varphi+\psi}$ gives
\begin{eqnarray*}|(f,g)_{\varphi + \frac{1-\ap}{2}\psi}|^2
&=& |(F_\psi\circ F_\psi^{-1}e^{\frac{1+\ap}{2}\psi}f,g)_{\varphi
+ \psi}|^2\\ &\leq&
(e^{\frac{1+\ap}{2}}f,e^{\frac{1+\ap}{2}}F_\psi^{-1}f)_{\varphi+\psi}(F_\psi
g,g)_{\varphi+\psi}\\ &\leq&
\frac{4(F_\psi^{-1}f,f)_{\varphi-\ap\psi}}{(1-\ap)^2}(\|e^{-\frac{1+\ap}{4}\psi}\dt_{\varphi+\frac{1-\ap}{2}\psi}g\|^2_{\varphi+\frac{1-\ap}{2}\psi}
+\|e^{-\frac{1+\ap}{4}\psi}dg\|^2_{\varphi+\frac{1-\ap}{2}\psi})
\end{eqnarray*}
where $F_\psi^{-1}$ is defined by \eqref{eq18}. Thus the estimate
\eqref{eq41} has been proved for $g\in {\rm Dom}(d^{*})\cap
C^\infty_p(\overline{\Omega})$. By using the density
lemma(proposition 1.2.4 in \cite{H1}), we know that \eqref{eq41}
holds for any $g\in {\rm Dom}(T^*)\cap{\rm Dom}(S).$ Consequently,
by the Lemma \ref{lem31}, there exists some $v\in
L_{p-1}^2(\Omega,\varphi+\frac{1-\ap}{2}\psi)$ such that
$$Tv=f,\ \ \
\|v\|_{\varphi+\frac{1-\ap}{2}\psi}^2\leq\frac{4}{(1-\ap)^2}(F_\psi^{-1}f,f)_{\varphi-\ap\psi}.
$$
Set $u=e^{-\frac{1+\ap}{4}\psi}v$, then we get
$u\in L_{p-1}^2(\Omega,\varphi-\ap\psi)$ and

\beq
du=f,\ \ \ \|u\|_{\varphi-\ap\psi}^2
=\|v\|_{\varphi+\frac{1-\ap}{2}\psi}^2
\leq\frac{4}{(1-\ap)^2}(F_\psi^{-1}f,f)_{\varphi-\ap\psi}. \label{eq44}
\eeq

Theorem \ref{thm41} now follows, in its full generality, from
\eqref{eq44} and the standard argument of smooth approximation
followed by taking weak limit as we did in the proof of \ref{thm31}.
\eprf

\bcor
\label{cor41}
Let $\Omega$ be a $p$-convex
domain in $\mbbr^n(1\leq p\leq n)$
and let $\varphi$ be a
$p$-plurisubharmonic function on $\Omega$, $\psi\in C^2(\Omega)$ be a
strictly plurisubharmonic function such that $-e^{-\psi}$ is $p$-plurisubharmonic.
For any constant $\ap\in [0,1)$ and
 closed $p$-form
$f\in L_p^2(\Omega,{\rm Loc})$, if
$$\int_\Omega\psi^{jk}f_{jK}
f_{kK}
e^{-\varphi+\ap\psi}<\infty$$
then there exists a $(p-1)$-form
$u\in L_{p-1}^2(\Omega,\varphi-\ap\psi)$
such that
$$ du=f,\ \
\|u\|_{\varphi-\ap\psi}^2\leq\frac{4}{p^2(1-\ap)^2}
\int_\Omega\psi^{jk}f_{jK} f_{kK}
e^{-\varphi+\ap\psi}$$
where $(\psi^{jk}):=(\psi_{jk})^{-1}$.
\ecor

\bprfc{cor41}  Corollary \ref{cor41} follows directly from Theorem
\ref{thm41} and the pointwise inequality \eqref{eq110}. \eprf

As a consequence of Theorem \ref{thm41}, we have the following analogue of  the Donnelly-Fefferman estimate(\cite{DF}).

\bthm \label{thm42} \emph{Let $\Omega$ be a $p$-convex domain in
$\mbbr^n(1\leq p\leq n)$ and let $\varphi$ be a $p$-plurisubharmonic
function on $\Omega$, $\psi\in C^2(\Omega)$ be a strictly
$p$-plurisubharmonic function such that $-e^{-\psi}$ is
$p$-plurisubharmonic. For any closed $p$-form $f\in
L_p^2(\Omega,{\rm Loc})$, if $$\int_\Omega\langle
F_\psi^{-1}f,f\rangle e^{-\varphi}<\infty$$ then there exists a
$(p-1)$-form $u\in L_{p-1}^2(\Omega,\varphi)$ such that
$$du=f,\ \
\|u\|_{\varphi}^2\leq 4\int_\Omega\langle F_\psi^{-1}f,f\rangle
e^{-\varphi}.$$} \ethm

\bprft{thm42}\ \ Theorem \ref{thm42} follows directly from Theorem
\ref{thm41} by choosing the constant $\ap$ to be $0$.\eprf

\bcor\label{cor42} Let $\Omega$ be a bounded $p$-convex domain in
$\mbbr^n(1\leq p\leq n)$ and let $\varphi$ be a $p$-plurisubharmonic
function on $\Omega$. For any closed $p$-form $f\in
L_{p}^2(\Omega,\varphi)$, there exists a $(p-1)$-form $u\in
L_{p-1}^2(\Omega,\varphi)$ such that
$$du=f,\ \
\|u\|_{\varphi}\leq \frac{2D}{p}\|f\|_{\varphi}.$$ where $D$ is the
diameter of $\Omega$. \ecor

\bprfc{cor42} \ \ Without loss of generality, we assume that
$\Omega$ contains the origin of $\mbbr^n$. Let
$$\psi=\frac{p|x|^2}{2D^2},$$ then \eqref{eq19} implies that
$$F_\psi^{-1}= \frac{D^2}{p^2}{\rm Id} \ {\rm  holds \ on} \ p{\rm -forms}.$$ Since the
Hessian of $-e^{-\psi}$ is given by
$$\frac{p}{D^2}e^{-\psi} (dx^i\otimes
dx^i-\frac{p}{D^2}x^idx^i\otimes x^jdx^j),$$ we know that any sum of
$p$ eigenvalues of the Hessian of $-e^{-\psi}$ is no less than
$$\frac{p}{D^2}e^{-\psi}[(1-\frac{p}{D^2}|x|^2)
+p-1]=\frac{p^2}{D^2}e^{-\psi}(1-\frac{|x|^2}{D^2})\geq 0.$$ So
$-e^{-\psi}$ is, by definition, a $p$-plurisubharmonic function on
$\Omega$(but not plurisubharmonic). Applying Theorem \ref{thm42}
with the weight function $\psi=\frac{p|x|^2}{2D^2},$ we obtain the following
estimate for the solution $u$
$$\|u\|_\varphi^2\leq\frac{4D^2}{p^2}\|f\|_\varphi^2.$$ This completes the proof of
Corollary \ref{cor42}.\eprf

\sec{Minimal $L^2$-solutions}

Let $\Omega$ be an open subset of $\mathbb{R}^n$ and $\varphi\in
L^{\infty}(\Omega,{\rm Loc})$, then the de Rham complex induces the
following complex of closed and densely defined
operators$$\cdots\rightarrow
L_{p-2}^2(\Omega,\varphi)\overset{d_{p-2}}{\rightarrow}
L_{p-1}^2(\Omega,\varphi)
\overset{d_{p-1}}{\rightarrow}L_{p}^2(\Omega,\varphi)\rightarrow\cdots,$$
where $d_{\ell}$'s denote the maximal(weak) differential operators
defined by the exterior derivatives. Then we have \beq \label{eq51}
{\rm Ker}d_{p-2,\varphi}^*\supseteq {\rm Ker}d_{p-1}^\perp \eeq and
since ${\rm Ker}d_{p-1}$ is a closed subspace of
$L_{p-1}^2(\Omega,\varphi)$ we also have the following orthogonal
decomposition

\beq \label{eq52}
L_{p-1}^2(\Omega,\varphi)={\rm
Ker}d_{p-1}^\perp\oplus{\rm Ker}d_{p-1}.
\eeq

Given a $d$-closed form $f\in L_{p}^2(\Omega,{\rm Loc})$, if there
is a $p$-form $u\in L_{p}^2(\Omega,\varphi)$ such that $du=f$, we
can decompose $u$ according to \eqref{eq52} \beq \label{eq53}
u=u_0+u_1\in ({\rm Dom}(d_{p-1})\cap{\rm
Ker}d_{p-1}^\perp)\oplus{\rm Ker}d_{p-1} \eeq which, together with
\eqref{eq51} above, implies that \beq \label{eq54} d_{p-1}u_0=f,\ \
d_{p-1,\varphi}^* u_0=0. \eeq

We will call the solution $u_0$ constructed in \eqref{eq53} the
\textbf{minimal solution} of $du=f$ in $L_{p-1}^2(\Omega,\varphi)$.\\

\brem\label{rem51} (i)For any $p$-convex open subset
$\Omega\subseteq\mathbb{R}^n$ and any closed $p$-form $f\in
L_p^2(\Omega,{\rm Loc})$, by Corollary \ref{cor31}, we can find some
$u\in L_{p-1}^2(\Omega,{\rm Loc})$ such that $du=f$. Let $\varphi\in
L^\infty(\Omega,{\rm Loc})$ and $\Omega^{'}\Subset\Omega$, previous
decomposition \eqref{eq53} applied to
$L_{p-1}^2(\Omega^{'},\varphi)$ gives the minimal solution of $du=f$
in $L_{p-1}^2(\Omega^{'},\varphi)$.\\ (ii)It is easy to see the
uniqueness of minimal solution, to be more precisely, by using
\eqref{eq53} we have $\|u_0\|_\varphi\leq\|u\|_\varphi$ holds for
any $u\in L_{p}^2(\Omega,\varphi)$ satisfying $du=f$, and the
equality holds if and only if $u=u_0$.\\ (iii)As an easy consequence
of (ii), we have the following monotonicity of $L^2$-solutions. Let
$\Omega_1\subseteq\Omega_2$ be open subsets of $\mathbb{R}^n$ and
$\varphi\in L^\infty(\Omega_2,{\rm Loc})$, for the minimal solution
$u_i$ of $du=f$ in $L_{p-1}^2(\Omega_i)(i=1,2)$, we have
$$\int_{\Omega_1}|u_1|^2e^{-\varphi}\leq
\int_{\Omega_2}|u_2|^2e^{-\varphi}.$$ Similarly, for any open set
$\Omega\subset\mathbb{R}^n, \varphi_i\in L^\infty(\Omega,{\rm Loc})$
and the minimal solution $u_i$ of $du=f$ in
$L_{p-1}^2(\Omega,\varphi_i)(i=1,2)$, if $\varphi_1\leq\varphi_2$ holds
on $\Omega$ then we have
$$\int_{\Omega}|u_1|^2e^{-\varphi_1}\ge
\int_{\Omega}|u_2|^2e^{-\varphi_2}.$$\erem

The minimal $L^2$-solution enjoys the following interior regularity
property.

\bpp\label{pp51} Under the conditions of Theorem \ref{thm31},  for
any $q\ge p$ and any closed $q$-form $f\in L^2_q(\Omega,\varphi)$
with $\int_{\Omega}\langle F_{\varphi}^{-1}f,f\rangle
e^{-\varphi}<\infty$, $du=f$ has a unique minimal solution $u_0$ in
$L^2_{q-1}(\Omega,\varphi)$, moreover if $f$ and the weight
$\varphi$ are both smooth then $u_0\in C_{q-1}^\infty(\Omega)$. The
same conclusion holds for Theorem \ref{thm41}.\epp \bprfp{pp51} The existence and uniqueness of minimal $L^2$-solution follows from the
decomposition \eqref{eq53} and Theorem \ref{thm31}. By \eqref{eq54},
we obtain
$$du_0=f,\ \ \delta_\varphi u_0=0$$ in the sense of distribution. This can be rewritten as $(d\delta_\varphi+\delta_\varphi d)u_0=\delta_\varphi f\in
C_{q-1}^\infty(\Omega)$. Now the smoothness of the minimal solution
$u_0$ follows from the interior elliptic regularity of the Hodge
Laplace operator $d\delta_\varphi+\delta_\varphi d$.  \eprf

If $\Omega$ is a strictly $p$-convex open set with smooth boundary,
it was proved(for compact Riemannian manifolds with smooth
$p$-convex boundary) in \cite{Sh} and \cite{W1} that
$\overline{\Omega}$ has the homotopy type of CW complex of dimension
$<p$. As an application of $L^2$-method we obtain the following
vanishing result of de Rham cohomology groups. Note that this result
was also obtained in \cite{AC}. We will generalize this result in section 8(see Proposition \ref{pp82} and Remark \ref{rem81}(ii)). 

\bcor \label{cor51}For any $p$-convex open subset
$\Omega\subseteq\mbbr^n(1\leq p\leq n)$, the de Rham cohomology
groups $H^q(\Omega,\mathbb{R})=0,p\leq q\leq n$. \ecor

\bprfc{cor51}Let $f\in C_q^\infty(\Omega)(p\leq q\leq n)$ be a
closed form. Since $p$-convexity implies $q$-convexity, there exists
a $q$-plurisubharmonic proper exhaustion function $\varphi$ such
that $\int_\Omega |f|e^{-\varphi}<\infty$. One can therefore find,
by Theorem \ref{thm31}(with the weight $\varphi(x)+|x|^2)$ and
proposition \ref{pp51}, a $(q-1)$-form $u\in C_{q-1}^\infty(\Omega)$
which solves the equation $du=f$ and this completes the proof of
$H^q(\Omega,\mathbb{R})=0$. \eprf

By the argument in \cite{Be2}(or \cite{Be4}), Prekopa's minimal principle follows from the estimate in Theorem \ref{thm31}(with $n=p=1$) applied to  the $L^2$-minimal solution given by  Proposition \ref{pp51} .
\bcor \label{cor52}  Let $\vphi(x,y)$ be a convex function in $\mbbr_x^{n}\times\mbbr_y^{m}$.
Define $\tilde{\vphi}$ by $$\tilde{\vphi}(x)=-\log\int_{\mbbr^{m}} e^{-\vphi(x,y)}dy.$$ Then $\tilde{\vphi}$ is a convex function on $\mbbr^{n}$.\ecor

We end this section by proving an estimate for $L^2$-minimal
solutions. The difference between this estimate and Theorem
\ref{thm41} is that the minimal solution here only depends on one of
the weights.  The
idea of the following proof goes back to \cite{Be3}and \cite{B1}(see also \cite{B4}, \cite{B5} and references therein).\\

\bthm \label{thm51} \emph{Let $\Omega$ be a $p$-convex domain in
$\mbbr^n(1\leq p\leq n)$ and let $\varphi$ be a $p$-plurisubharmonic
function on $\Omega$, $\psi\in C^2(\Omega)$. If we assume, in
addition, that there are a function $0\leq\omega<1$ and a constant
$\alpha\in [0,1)$ such that the quadratic form $\omega^2
D^2\psi-d\psi\otimes d\psi$ is $p$-positive semi-definite on
$\Omega$ and that $\omega\leq\alpha$ holds on supp$f$, where $f\in
L_{p}^2(\Omega,{\rm Loc})$ is a closed $p$-form, then the minimal
solution, denoted by $u_\varphi$, of $du=f$ in
$L_{p-1}^2(\Omega,\varphi)$ satisfies
$$\int_\Omega(1-\omega^2)|u_\varphi|^2e^{-\varphi+\psi}\leq\frac{1+\alpha}{1-\alpha}\int_\Omega\langle
F_\psi^{-1}f,f\rangle e^{-\varphi+\psi}$$where
$D^2\psi:=\psi_{ij}dx^i\otimes dx^j$ is the hessian of $\psi$.}
\ethm

\bprft{thm51} By the monotonicity discussed in remark \ref{rem51}
(iii) and the standard argument of approximation followed by taking
weak limit, we can assume in addition  that $\Omega$ is a bounded
open set with smooth boundary and that $\varphi,\psi$ are both smooth up to
the boundary of $\Omega$. Set $$u=e^{\psi}u_\varphi,$$ by
\eqref{eq53}, $u$ is the minimal solution of
$du=e^{\psi}(d\psi\wedge u_\varphi+f):=e^{\psi}g$ in
$L_{p-1}^2(\Omega,\varphi+\psi)$. Since the quadratic form $\omega^2
D^2\psi-d\psi\otimes d\psi$ is $p$-positive semi-definite and
$\omega\leq\alpha$ on supp$f$, by using Lemma \ref{lem11} to
$$\theta=\omega^2 D^2\psi \ {\rm and}  \ \tau=d\psi,$$ it follows that
$F_\psi^{-1}(d\psi\wedge u_\varphi)$ is well-defined and \beq
\langle F_{\varphi+\psi}^{-1} g,g\rangle &\leq & \langle
F_{\psi}^{-1} d\psi\wedge u_\varphi, d\psi\wedge
u_\varphi\rangle+2\langle F_\psi^{-1} f, d\psi\wedge
u_\varphi\rangle+\langle F_\psi^{-1} f,
f\rangle \nno \\
&\leq&\omega^2|u_\varphi|^2+2\langle f,F_\psi^{-1}
f\rangle^{\frac{1}{2}}\cdot\alpha|u_\varphi|+\langle
F_\psi^{-1} f, f\rangle \nno \\
&\leq
&\frac{\alpha+\omega^2}{1+\alpha}|u_\varphi|^2+\frac{1}{1-\alpha}\langle
F_\psi^{-1} f, f\rangle.\label{eq55} \eeq
Since $\varphi+\psi$ is $p$-plurisubharmonic and  $du=e^\psi g$, we can apply Theorem
\ref{thm31} to get
\beq\int_\Omega|u_\varphi|^2e^{-\varphi+\psi}&=&\|u\|_{\varphi+\psi}^2\nno
\\&\leq& \int_\Omega\langle F_{\varphi+\psi}^{-1}(e^\psi g),e^\psi g\rangle
e^{-\varphi-\psi}\nno \\&=& \int_\Omega\langle F_{\varphi+\psi}^{-1}g,g\rangle
e^{-\varphi+\psi}\nno\\ &\leq&
\int_\Omega \frac{\alpha+\omega^2}{1+\alpha}|u_\varphi|^2e^{-\varphi+\psi}+\frac{1}{1-\alpha}\int_\Omega\langle
F_\psi^{-1} f, f\rangle e^{-\varphi+\psi}\label{eq56}\eeq where we
have used the inequality \eqref{eq55}. Now the desired
$L^2$-estimate
follows directly from \eqref{eq56}. \\
\eprf

\brem\label{rem52} Theorem \ref{thm51} could be used to deduce a
weaker version of Theorem \ref{thm41}. Let $\varphi$ be a
$p$-plurisubharmonic function on $\Omega$, $\psi\in C^2(\Omega)$ be
a function such that $-e^{-\psi}$ is $p$-plurisubharmonic. Then for
any constant $\ap\in [0,1)$, $\alpha\psi$ satisfies the conditions
assumed in Theorem \ref{thm51} with $\omega$ given by the constant
$\sqrt{\alpha}$, and consequently we obtain
$$\|u\|_{\varphi-\ap\psi}^2\leq\frac{1}{\alpha(1-\sqrt{\ap})^2}\int_\Omega\langle
F_\psi^{-1}f,f\rangle e^{-\varphi+\ap\psi}.$$\erem

\sec{Non-plurisubharmonic weights}

Next we prove a theorem which  has the feature of allowing
non-plurisubharmonic weights. This kind of result will provide more
flexibility in choosing weights for $L^2$-estimates. Such an estimate for $\bar\partial$-problem was proved by B{\L}ocki(\cite{B4},\cite{B5}).

\bthm\label{thm61} Let $\Omega$ be a $p$-convex domain in
$\mathbb{R}^n(1\leq p\leq n)$ and let $\varphi\in C^2(\Omega)$ be a
$p$-plurisubharmonic function on $\Omega$ and $\psi\in C^1(\Omega).$
There are a function $0\leq\omega<2$ and a constant $\alpha\in
[0,2)$ such that the quadratic form $\omega^2
D^2\varphi-d\psi\otimes d\psi$ is $p$-positive semi-definite on
$\Omega$ and that $\omega\leq\alpha$ on supp$f$ where $f\in
L_{p}^2(\Omega,{\rm Loc})$ is a closed $p$-form. If
$$\int_\Omega\langle F_\varphi^{-1}f,f\rangle
e^{-\varphi+\psi}<\infty,$$then there exists a $(p-1)$-form $u\in
L_{p-1}^2(\Omega,\varphi-\psi)$ such that
$$du=f,\ \
\int_\Omega(1-\frac{\omega^2}{4})|u|^2
e^{-\varphi+\psi}\leq\frac{2+\alpha}{2-\alpha}\int_\Omega \langle
F_\varphi^{-1} f, f\rangle e^{-\varphi+\psi},$$where
$F_\varphi^{-1}$ is defined by \eqref{eq18} and it is required
implicitly that $F_\varphi^{-1}f$ is defined almost everywhere in
$\Omega$, $D^2\varphi:=\varphi_{ij}dx^i\otimes dx^j$ is the Hessian
of $\varphi$. \ethm

\bprft{thm61}By the standard argument used in the proof of Theorem
\ref{thm31}, we may assume, without loss of generality, that
$\Omega$ is a bounded open set with smooth $p$-convex boundary and
that $\varphi, \psi$ are both smooth up to the boundary. In this
case, there exists a unique minimal solution, denoted by $u_0$, of
$du=f$ in $L_{p-1}^2(\Omega,\varphi-\frac{1}{2}\psi)$. For $u_0$, we
have \beq \label{eq61} \int_\Omega\langle u_0, v\rangle
e^{-\varphi+\frac{1}{2}\psi}=0 \eeq for any closed $(p-1)$-form
$v\in L_{p-1}^2(\Omega,\varphi-\frac{1}{2}\psi)$. Set
$$u=e^{\frac{1}{2}\psi}u_0$$ then \eqref{eq61} implies that $u$ is the
minimal solution of $du=g$ in $L_{p-1}^2(\Omega,\varphi)$ where $g$
is the closed $p$-form given by \beq
g=e^{\frac{1}{2}\psi}(\frac{1}{2}d\psi\wedge u_0+f).\nno \eeq

By lemma \ref{lem11}, $F_\varphi g$ is well-defined and we have the
following pointwise inequality \beq \langle F_\varphi^{-1} g,
g\rangle&=&(\frac{1}{4}\langle F_\varphi^{-1} d\psi\wedge u_0,
d\psi\wedge u_0\rangle+\langle F_\varphi^{-1} f, d\psi\wedge
u_0\rangle+\langle F_\varphi^{-1} f,
f\rangle)e^\psi \nno \\
&\leq&(\frac{\omega^2}{4}|u_0|^2+\langle f,F_\varphi^{-1}
f\rangle^{\frac{1}{2}}\cdot\alpha|u_0|+\langle F_\varphi^{-1} f,
f\rangle)e^\psi  \label{eq62} \eeq where we have used the
assumptions that $\omega^2 D^2\varphi-d\psi\otimes d\psi$ is
$p$-positive semi-definite and that $\omega\leq\alpha$ holds on
${\rm supp}f$.

Since $\varphi$ is by assumption a $p$-plurisubharmonic function,
from Theorem \ref{thm31} it follows
that$$\int_\Omega|u_0|^2e^{-\varphi+\psi}=\|u\|_\varphi^2\leq\int_\Omega\langle
F_\varphi^{-1} g, g\rangle e^{-\varphi},$$which, together with
\eqref{eq62}, implies
\begin{eqnarray*} \int_\Omega(1-\frac{\omega^2}{4})|u_0|^2e^{-\varphi+\psi}&\leq&\int_\Omega (\langle
f,F_\varphi^{-1} f\rangle^{\frac{1}{2}}\cdot\alpha|u_0|+\langle
F_\varphi^{-1} f, f\rangle)e^{-\varphi+\psi}\\
&\leq&\epsilon\int_\Omega(1-\frac{\omega^2}{4})|u_0|^2e^{-\varphi+\psi}
+\int_\Omega [1+\frac{\alpha^2}{(4-\omega^2)\epsilon}]\langle
F_\varphi^{-1} f, f\rangle
e^{-\varphi+\psi}\\&\leq&\epsilon\int_\Omega(1-\frac{\omega^2}{4})|u_0|^2e^{-\varphi+\psi}
+\int_\Omega [1+\frac{\alpha^2}{(4-\alpha^2)\epsilon}]\langle
F_\varphi^{-1} f, f\rangle e^{-\varphi+\psi}
\end{eqnarray*}where $0<\epsilon<1$ is any constant. Set
$$\epsilon=\frac{\alpha}{2+\alpha},$$ the above inequality gives
$$\int_\Omega(1-\frac{\omega^2}{4})|u_0|^2
e^{-\varphi+\psi}\leq\frac{2+\alpha}{2-\alpha}\int_\Omega \langle
F_\varphi^{-1} f, f\rangle e^{-\varphi+\psi},$$ hence $u_0$ is the
desired solution. \eprf

As an immediate consequence of the above theorem, if the function
$\omega$ is constant we have the following corollary\\

\bcor \label{cor61}Let $\Omega$ be a $p$-convex domain in
$\mathbb{R}^n(1\leq p\leq n)$ and let $\varphi\in C^2(\Omega)$ be a
$p$-plurisubharmonic function on $\Omega$ and $\psi\in C^1(\Omega).$
There is a constant $\alpha\in [0,2)$ such that the symmetric
bilinear form $\alpha^2D^2\varphi-d\psi\otimes d\psi$ is
$p$-positive semi-definite on $\Omega$. For any closed $p$-form
$f\in L_{p}^2(\Omega,{\rm Loc})$, if
$$\int_\Omega\langle F_\varphi^{-1}f,f\rangle
e^{-\varphi+\psi}<\infty,$$then there exists a $(p-1)$-form $u\in
L_{p-1}^2(\Omega,\varphi-\psi)$ such that
$$du=f,\ \
\|u\|_{\varphi-\psi}^2\leq\frac{4}{(2-\alpha)^2}\int_\Omega\langle
F_\varphi^{-1}f,f\rangle e^{-\varphi+\psi}.$$
\ecor

\brem\label{rem61}(i)If we choose the constant $\alpha=0$ and the
weight function $\psi=0$, then Corollary \ref{cor61} recovers Theorem \ref{thm31}.\\
(ii) We can give an alternative proof of Theorem \ref{thm41} by
using Corollary \ref{cor61} in the following way. Let
$\varphi_1=\varphi+\psi$ and $\psi_1=(1+\alpha)\psi$, then
$\varphi_1$ is $p$-plurisubharmonic. Since
$$(1+\alpha)^2D^2\varphi_1-
d\psi_1\otimes d\psi_1=(1+\alpha)^2[D^2\varphi+
e^{\psi}D^2(-e^{-\psi})],$$ the assumption that $\varphi$ and
$-e^{-\psi}$ are both $p$-plurisubharmonic functions implies that
$(1+\alpha)^2D^2\varphi_1- d\psi_1\otimes d\psi_1$ is $p$-positive
semi-definite. Applying Corollary \ref{cor61} to the
weights $\varphi_1$ and $\psi_1$, we obtain Theorem \ref{thm41}.\\
(iii)The proof of Theorem \ref{thm41} given in (ii) does not
indicate the estimate $(*)$ in section 4. Actually, Corollary 6.1
also follows from the following estimate whose proof is an imitation
of that of $(*)$. Let $\varphi\in C^\infty(\overline{\Omega})$ be a
$p$-plurisubharmonic function and let $\psi\in
C^\infty(\overline{\Omega}))$ be a function such that the symmetric
form $\alpha D^2\varphi-d\psi\otimes d\psi$ is $p$-positive
semi-definite for some constant $\alpha\in [0,2)$, we have the
following apriori estimate
$$\|\delta_{\varphi-\frac{1}{2}\psi}g\|^2_{\varphi}
+\|dg\|^2_{\varphi}\geq \frac{(2-\alpha)^2}{4}\int_\Omega\langle
F_{\varphi} g,g\rangle e^{-\varphi},\eqno{(**)}$$ for any $p$-form
$g\in {\rm Dom}(d^*)\cap C^\infty_{p}(\overline{\Omega})$ on
$p$-convex domains with smooth boundary.\erem

\sec{$L^2$-estimates on $p$-convex Riemannian manifolds}

We will generalize  the results established in sections 2-6 to
Riemannian manifolds. To this end, we only need to take care of the
curvature term which enters the apriori estimate and we will focus
on such modifications.

Let $(M,ds^2)$ be an oriented Riemannian manifolds of dimension $n$.
We denote by $R_{XY}=D_XD_Y-D_YD_X-D_{[X,Y]}$ the curvature of the
Levi-Civita connection $D$.  Let $\{e_1,\cdots,e_n\}$ be locally
defined orthonormal frame field of the tangent bundle and
$\{\omega^1,\cdots,\omega^n\}$ be its dual coframe field. Since $D$
is torsion free, the exterior differential operator $d$ and its
formal adjoint $\delta$ satisfy \beq \label{eq71}d=\omega^i\wedge
D_{e_i}, \ \ \delta=-e_i\lrcorner D_{e_i}.\eeq For any $\varphi\in
C^\infty(M)$, we denote as before
$$F_\varphi=\varphi_{ij}\omega^j\wedge e_i\lrcorner$$where
$\varphi_{ij}$'s are given by the Hessian
$D^2\varphi:=\varphi_{ij}\omega^i\otimes\omega^j$ of $\varphi$.

For our later use, we collect here some easy geometric computations.

\blem \label{lem71} For any $\varphi\in C^\infty(M)$ and any
$p$-form $g\in C_p^\infty(M)$, we have the following
identity\beq\label{eq72}
L_{\nabla\varphi}g=D_{\nabla\varphi}g+F_\varphi g\eeq where
$L_{\nabla\varphi}=d\nabla\varphi\lrcorner+\nabla\varphi\lrcorner d$
is the Lie derivative and $\nabla\varphi$ is the gradient of
$\varphi$. \elem

\bprfl{lem71}By repeated use of the first formula in \eqref{eq71},
we have \begin{eqnarray*} \nabla\varphi\lrcorner dg
&=&\nabla\varphi\lrcorner \omega^i\wedge D_{e_i}g \\ &=&
\langle\nabla\varphi,e_i\rangle
D_{e_i}g-\omega^i\wedge\nabla\varphi\lrcorner D_{e_i}g
\\&=&D_{\nabla\varphi}g-\omega^i\wedge[D_{e_i}(\nabla\varphi\lrcorner g)-(D_{e_i}\nabla\varphi)\lrcorner
g]\\
&=&D_{\nabla\varphi}g-\omega^i\wedge D_{e_i}(\nabla\varphi\lrcorner
g)+\varphi_{ij}\omega^i\wedge e_j\lrcorner g\\ &=&
D_{\nabla\varphi}g + F_\varphi g-d\nabla\varphi\lrcorner
g.\end{eqnarray*} The proof is complete. \eprf

\blem \label{lem72} Let $\Omega\Subset M$ be an open subset with
smooth boundary. For any differential forms $f\in
C^\infty_{p+1}(\overline{\Omega}),g\in
C^\infty_{p}(\overline{\Omega})$, we have the following identities
\beq\label{eq73} \int_\Omega\langle f,dg\rangle&=&\int_\Omega\langle
\delta f,g\rangle + \int_{\partial\Omega}\langle\nabla\rho\lrcorner
f,g\rangle\frac{1}{|\nabla\rho|}\\ \label{eq74}
\int_{\Omega}\langle\triangle
g,g\rangle&=&-\int_\Omega|Dg|^2+\int_{\partial\Omega}\langle
D_{\nabla\rho} g,g\rangle\frac{1}{|\nabla\rho|}\eeq where
$\triangle:={\rm tr}D^2$ is the Laplacian, $\rho$ is a defining
function for $\Omega$, i.e., $\rho\in C^\infty(\overline{\Omega})$
satisfying $\rho<0$ in $\Omega$, $\rho=0$ and $\nabla\rho\neq 0$ on
$\partial\Omega$. \elem \bprfl{lem72} Set $$X=\langle g,e_i\lrcorner
f\rangle e_i,\ \ Y=\langle g,D_{e_i} g\rangle e_i,$$ it is obvious
that $X,Y$ are both well-defined smooth vector fields on
$\overline{\Omega}$. By using \eqref{eq71}, we see that
\beq\label{eq75} {\rm div}X=\langle dg,f\rangle-\langle g,\delta
f\rangle\eeq and that \beq\label{eq76}{\rm div}Y=\langle \triangle
g,g\rangle+|Dg|^2.\eeq Now the divergence theorem gives the required
identities \eqref{eq73}, \eqref{eq74} by integrating \eqref{eq75}
and \eqref{eq76} respectively.\eprf

We use the same notation $d$ to denote the maximal(weak)
differential operator $d: L_{p-1}^2(\Omega)\rightarrow
L_p^2(\Omega)$ where $\Omega\Subset M$ be an smooth open subset. We
also denote the adjoint of the closed and densely operator by $d^*$.
Form \eqref{eq73}, it is easy to see that \beq\label{eq77}
C_p^\infty(\overline{\Omega})\cap{\rm Dom}(d^*)=\{g\in
C_p^\infty(\overline{\Omega}) \mid \nabla\rho\lrcorner g=0\ {\rm on}
\ \partial\Omega\}\eeq where $1\leq p\leq n$.\\

To establish the basic estimate in section 3 on Riemannian
manifolds, we first compute the integral $\int_\Omega |dg|^2+|\delta
g|^2$ for any $g\in C_p^\infty(\overline{\Omega})\cap{\rm
Dom}(d^*)$. From \eqref{eq77} and lemma \ref{lem72}, it follows that
\beq\label{eq78}\int_\Omega |dg|^2+|\delta g|^2&=&\int_\Omega
\langle(d\delta+\delta d)g,g\rangle+\int_{\partial\Omega}
\Big(\langle\nabla\rho\lrcorner dg,g\rangle-\langle\nabla\rho\lrcorner
g,\delta g\rangle\Big)\frac{1}{|\nabla\rho|}\nno
\\&=&\int_\Omega
\langle(d\delta+\delta d)g,g\rangle+\int_{\partial\Omega}
\langle\nabla\rho\lrcorner dg,g\rangle\frac{1}{|\nabla\rho|}.\eeq
Let us choose the orthonormal frame field $\{e_1,\cdots,e_n\}$ to be
adapted to $\partial\Omega$ with
$$e_n=\frac{\nabla\rho}{|\nabla\rho|},$$ then we know by \eqref{eq77}
that \beq\label{eq79} \langle d\nabla\rho\lrcorner g,g\rangle =
\sum_{\nu=1}^{n-1}\langle\omega^\nu\wedge
D_{e_\nu}(\nabla\rho\lrcorner g),g\rangle+\langle
D_{e_n}(\nabla\rho\lrcorner g),e_n\lrcorner g\rangle=0\eeq holds on
the boundary $\partial\Omega$. Combining \eqref{eq72},\eqref{eq78}
and \eqref{eq79} gives the next identity
\beq\label{eq710}\int_\Omega |dg|^2+|\delta g|^2&=&\int_\Omega
\langle(d\delta+\delta d)g,g\rangle+\int_{\partial\Omega} \langle
L_{\nabla\rho}g,g\rangle\frac{1}{|\nabla\rho|}\nno \\ &=&
\int_\Omega \langle(d\delta+\delta d)g,g\rangle+
\int_{\partial\Omega} \langle D_{\nabla\rho}g+F_\rho
g,g\rangle\frac{1}{|\nabla\rho|}. \eeq To handle the first term on
the right hand side of \eqref{eq710}, we use the
Bochner-Weitzenb\"{o}ck formula \beq\label{eq711} (d\delta+\delta
d)g=-\triangle g+\omega^j\wedge e_i\lrcorner R_{e_ie_j}g.\eeq Recall
that the curvature operator
$\mathfrak{R}:\bigwedge^2T^*M\rightarrow\bigwedge^2T^*M$ is defined
as a self-adjoint linear map by
$$\mathfrak{R}(\omega^i\wedge\omega^j):=R_{ij\ell
k}\omega^k\wedge\omega^\ell$$where $R_{ijk\ell}:=\langle
R_{e_ie_j}e_k,e_\ell\rangle.$ It is known that(cf.
\cite{W2})\beq\label{eq712}\langle\omega^j\wedge e_i\lrcorner
R_{e_ie_j}g,g\rangle=\sum_{i_1<\cdots<i_p}\langle\mathfrak{R}\xi^g_{i_1\cdots
i_p},\xi^g_{i_1\cdots i_p}\rangle\eeq where $\xi^g_{i_1\cdots i_p}$
is the $2$-form given by \beq\label{eq713}\xi^g_{i_1\cdots
i_p}=\sum_{a=1}^p\sum_{i=1}^ng_{i_1\cdots (i)_a\cdots
i_p}\omega^i\wedge\omega^{i_a}.\eeq By \eqref{eq74},\eqref{eq711}
and \eqref{eq712}, it is easy to see the following
equality\beq\label{eq714}\int_\Omega \langle(d\delta+\delta
d)g,g\rangle=\int_\Omega
|Dg|^2+\sum_{i_1<\cdots<i_p}\langle\mathfrak{R}\xi^g_{i_1\cdots
i_p},\xi^g_{i_1\cdots i_p}\rangle-\int_{\partial\Omega}\langle
D_{\nabla\rho}g,g\rangle\frac{1}{|\nabla\rho|}.\eeq Substituting
\eqref{eq714} into \eqref{eq710} implies that
\beq\label{eq715}\int_\Omega |dg|^2+|\delta g|^2=\int_\Omega
|Dg|^2+\sum_{i_1<\cdots<i_p}\langle\mathfrak{R}\xi^g_{i_1\cdots
i_p},\xi^g_{i_1\cdots i_p}\rangle+ \int_{\partial\Omega} \langle
F_\rho g,g\rangle\frac{1}{|\nabla\rho|}\eeq holds for any $g\in
C^\infty_{p}(\overline{\Omega})\cap {\rm Dom}(d^*)$.\\

Now we commence introducing a weight function into the identity
\eqref{eq715}.

\blem \label{lem73}Let $\Omega\Subset M$ be an open subset with a
defining function $\rho\in C^\infty(\overline{\Omega})$ ,
$\varphi\in C^\infty(\overline{\Omega})$. Then for any differential
form $g\in C^\infty_{p}(\overline{\Omega})\cap {\rm Dom}(d^*)$, we
have the following identities \beq\label{eq716}\int_\Omega
\Big(|dg|^2+|\delta_\varphi g|^2\Big)e^{-\varphi}&=&\int_\Omega
\Big(|Dg|^2+\sum_{i_1<\cdots<i_p}\langle\mathfrak{R}\xi^g_{i_1\cdots
i_p},\xi^g_{i_1\cdots i_p}\rangle+\langle F_\varphi
g,g\rangle\Big)e^{-\varphi}\nno \\&& +\int_{\partial\Omega} \langle
F_\rho g,g\rangle\frac{e^{-\varphi}}{|\nabla\rho|}\eeq where
$\delta_\varphi:=e^\varphi\circ\delta\circ e^{-\varphi}$ is the
formal adjoint of $d$ with respect to the weight $\varphi$ and
$\xi^g_{i_1\cdots i_p}$ is defined by \eqref{eq713}.\elem
\bprfl{lem73}For any $g\in C^\infty_{p}(\overline{\Omega})\cap {\rm
Dom}(d^*)$, set
$$h=e^{-\frac{\varphi}{2}}g$$ then we know by \eqref{eq77} $h\in
C^\infty_{p}(\overline{\Omega})\cap {\rm Dom}(d^*)$. The equality
\eqref{eq715} applied to $h$ gives \beq\label{eq717} \int_\Omega
(|dg|^2+|\delta_\varphi g|^2)e^{-\varphi}&=&
\int_\Omega|dh+\frac{1}{2}d\varphi\wedge h|^2+|\delta
h+\frac{1}{2}\nabla\varphi\lrcorner h|^2\nno \\ &=&
\int_\Omega|dh|^2+|\delta h|^2+\langle
L_{\nabla\varphi}h,h\rangle\nno \\ &&+\frac{1}{4}(|d\varphi\wedge
h|^2+|\nabla\varphi\lrcorner h|^2)\nno \\ &=& \int_\Omega
|Dh|^2+\sum_{i_1<\cdots<i_p}\langle\mathfrak{R}\xi^h_{i_1\cdots
i_p},\xi^h_{i_1\cdots i_p}\rangle+\langle F_\varphi h,h\rangle\nno
\\&&+ \int_\Omega \langle D_{\nabla\varphi}h,h\rangle+
\frac{1}{4}|d\varphi|^2|h|^2\nno \\&&+\int_{\partial\Omega} \langle F_\rho
h,h\rangle\frac{1}{|\nabla\rho|}\eeq where we have also used
\eqref{eq73} to get the second equality, \eqref{eq72} and the
Lagrange identity to get the last equality. By substituting
$h=e^{\frac{-\varphi}{2}}g$ into \eqref{eq717}, we obtain the
desired identity as follows \beq \int_\Omega \Big(|dg|^2+|\delta_\varphi
g|^2\Big)e^{-\varphi}&=& \int_\Omega \Big(|Dg-\frac{1}{2}d\varphi\otimes
g|^2+\sum_{i_1<\cdots<i_p}\langle\mathfrak{R}\xi^g_{i_1\cdots
i_p},\xi^g_{i_1\cdots i_p}\rangle\Big)e^{-\varphi}\nno \\
&&+\int_\Omega\Big(\langle F_\varphi g,g\rangle+\langle
D_{\nabla\varphi}g-\frac{1}{2}|\nabla\varphi|^2g,g\rangle\Big)e^{-\varphi}\nno
\\ &&+\frac{1}{4}\int_\Omega|d\varphi|^2|g|^2e^{-\varphi}+\int_{\partial\Omega} \langle F_\rho
g,g\rangle\frac{e^{-\varphi}}{|\nabla\rho|}\nno \\ &=& \int_\Omega
\Big(|Dg|^2+\sum_{i_1<\cdots<i_p}\langle\mathfrak{R}\xi^g_{i_1\cdots
i_p},\xi^g_{i_1\cdots i_p}\rangle+\langle F_\varphi g,g\rangle\Big)e^{-\varphi}\nno \\
&&+\int_\Omega\Big(\langle
D_{\nabla\varphi}g-\frac{1}{2}|\nabla\varphi|^2g,g\rangle-\langle
Dg,d\varphi\otimes g\rangle\Big)e^{-\varphi}\nno
\\ &&+\frac{1}{2}\int_\Omega|d\varphi|^2|g|^2e^{-\varphi}+\int_{\partial\Omega} \langle F_\rho
g,g\rangle\frac{e^{-\varphi}}{|\nabla\rho|}\nno \\ &=& \int_\Omega
\Big(|Dg|^2+\sum_{i_1<\cdots<i_p}\langle\mathfrak{R}\xi^g_{i_1\cdots
i_p},\xi^g_{i_1\cdots i_p}\rangle+\langle F_\varphi
g,g\rangle\Big)e^{-\varphi}\nno
\\ &&+\int_{\partial\Omega} \langle F_\rho
g,g\rangle\frac{e^{-\varphi}}{|\nabla\rho|}.\nno\eeq The proof is
complete. \eprf

Before we prove the $L^2$-existence theorem on $(M,ds^2)$, we need
to bound the curvature term in \eqref{eq716}. Set \beq\label{eq718}
\lambda_{\mathfrak{R}}(x)&:=&{\rm the \ smallest \
eigenvalue \ of} \ \mathfrak{R}(x)\nno \\
\Lambda_{\mathfrak{R}}(x)&:=&{\rm the \ largest \ eigenvalue \ of} \
\mathfrak{R}(x)\eeq for any $x\in M$. Then we have, for any $p$-form
$g$, the following pointwise inequalities for the curvature term
$\sum_{i_1<\cdots<i_p}\langle\mathfrak{R}\xi^g_{i_1\cdots
i_p},\xi^g_{i_1\cdots i_p}\rangle$ in \eqref{eq716}.

\blem \label{lem74} \beq\label{eq719}
p(n-p)\lambda_{\mathfrak{R}}|g|^2\leq\sum_{i_1<\cdots<i_p}\langle\mathfrak{R}\xi^g_{i_1\cdots
i_p},\xi^g_{i_1\cdots i_p}\rangle\leq
p(n-p)\Lambda_{\mathfrak{R}}|g|^2\eeq where the $\xi^g_{i_1\cdots
i_p}$'s are defined by \eqref{eq713}.\elem \bprfl{lem74} By
definition of $\xi^g_{i_1\cdots i_p}$, we get
\beq\sum_{i_1<\cdots<i_p}\langle\mathfrak{R}\xi^g_{i_1\cdots
i_p},\xi^g_{i_1\cdots i_p}\rangle&\geq&
\lambda_{\mathfrak{R}}\sum_{i_1<\cdots<i_p}[\sum_{a=1}^p\sum_{i=1}^ng^2_{i_1\cdots
(i)_a\cdots i_p}-\sum_{a,b=1}^pg_{i_1\cdots (i_b)_a\cdots
i_p}g_{i_1\cdots (i_a)_b\cdots i_p}]\nno \\ &=&
\lambda_{\mathfrak{R}}\sum_{i_1<\cdots<i_p}[\sum_{a=1}^p\sum_{i=1}^ng^2_{i_1\cdots
(i)_a\cdots i_p}-pg^2_{i_1\cdots i_p}]\nno \\ &=&
\lambda_{\mathfrak{R}}\sum_{i_1<\cdots<i_p}\sum_{a=1}^p\sum_{i\neq
i_1,\cdots,i_p}g^2_{i_1\cdots (i)_a\cdots i_p}\nno
\\ &=&\lambda_{\mathfrak{R}}\sum_{j_1<\cdots<j_p}\sum_{i_1<\cdots<i_p}\sum_{a=1}^p\sum_{i\neq
i_1,\cdots,i_p}{\rm sgn}\left(\begin{array}{c}
i_1\cdots(i)_a\cdots i_p\\
j_{1}\cdots\cdots\cdots j_{p}\end{array}\right)^2g^2_{j_1\cdots
j_p}.\nno\eeq where sgn denotes the signature of permutation and we
have used the identity \beq g_{i_1\cdots (i)_a\cdots
i_p}=\sum_{j_1<\cdots<j_p}{\rm sgn}\left(\begin{array}{c}
i_1\cdots(i)_a\cdots i_p\\
j_{1}\cdots\cdots\cdots j_{p}\end{array}\right)g_{j_1\cdots
j_p}.\nno \eeq For fixed $j_1<\cdots<j_p$ we have \beq
\sum_{i_1<\cdots<i_p}\sum_{a=1}^p\sum_{i\neq i_1,\cdots,i_p}{\rm
sgn}\left(\begin{array}{c}
i_1\cdots(i)_a\cdots i_p\\
j_{1}\cdots\cdots\cdots j_{p}\end{array}\right)^2=p(n-p),\nno \eeq
because only terms given by
$\{i_1,\cdots,i_p\}=\{j_1,\cdots,\widehat{j_a},\cdots,j_p,k\}(k\neq
j_1,\cdots,j_p \ {\rm and} \ 1\leq a\leq p)$ contribute to the sum.
We can therefore rewrite the above inequality as
\beq\sum_{i_1<\cdots<i_p}\langle\mathfrak{R}\xi^g_{i_1\cdots
i_p},\xi^g_{i_1\cdots i_p}\rangle\geq
p(n-p)\lambda_\mathfrak{R}\sum_{j_1<\cdots<j_p}g^2_{j_1\cdots
j_p}.\nno \eeq The second inequality can be proved in the same
manner, and the proof is complete. \eprf

In order to establish $L^2$-existence theorem, we will use
$d:L^2_{p-1}(\Omega,\varphi)\rightarrow L^2_{p}(\Omega,\varphi)$,
the maximal(weak) differential operator between the weighted
$L^2$-spaces. Let $d_\varphi^*:L^2_{p}(\Omega,\varphi)\rightarrow
L^2_{p-1}(\Omega,\varphi)$ be the adjoint operator. As mentioned
before, we know by \eqref{eq73} that the formal adjoint of $d$ w.r.t
the weight is given by $\delta_\varphi$, and consequently we have
\beq\label{eq720}C^\infty_p(\overline{\Omega})\cap{\rm
Dom}(d_\varphi^*)=C^\infty_p(\overline{\Omega})\cap{\rm
Dom}(d^*)=\{g\in C_p^\infty(\overline{\Omega}) \mid
\nabla\rho\lrcorner g=0\ {\rm on} \ \partial\Omega\}.\eeq

Now we are in the position to prove the main result of this section.

\bthm \label{thm71}Let $(M,ds^2)$ be a $n$-dimensional oriented
$p$-convex Riemannian manifold. Let $\varphi\in C^2(M)$ be a
$p$-plurisubharmonic function on $M$. If
$F_\varphi+p(n-p)\lambda_\mathfrak{R}{\rm Id}$ is $p$-positive
semi-definite on $M$, then for any closed $p$-form $f\in
L^2_p(M,{\rm Loc})$ with
$$\int_M\langle[F_\varphi+p(n-p)\lambda_\mathfrak{R}{\rm
Id}]^{-1}f,f\rangle e^{-\varphi}<\infty,$$there exists some
$(p-1)$-form $u \in L^2_{p-1}(M,\varphi)$ such that $$du=f  \ {\rm
and} \ \int_M|u|^2e^{-\varphi}\leq
\int_M\langle[F_\varphi+p(n-p)\lambda_\mathfrak{R}{\rm
Id}]^{-1}f,f\rangle e^{-\varphi}$$where $1\leq p\leq n$,
$[F_\varphi+p(n-p)\lambda_\mathfrak{R}{\rm Id}]^{-1}$ is defined by
\eqref{eq18} and $\lambda_\mathfrak{R}$ is given by \eqref{eq718}.
Moreover, if $f$ and $\varphi$ are both assumed additionally to be
smooth then we can choose $u$ to be a smooth form. \ethm
\bprft{thm71} It has been proved in \cite{HL3} that $M$ admits a
smooth $p$-plurisubharmonic proper exhaustion function, so $M$
itself can be exhausted by  compact open sunsets with smooth
$p$-convex boundary. Since the resulting $L^2$-estimate enables us
to apply the standard argument of approximation to take weak limit,
we only need to work on a smooth domain $\Omega\Subset M$ which has
$p$-convex boundary. From \eqref{eq716}, \eqref{eq719} and
\eqref{eq720}, it follows that
\beq\label{eq721}\int_\Omega(|dg|^2+|\delta_\varphi
g|^2)e^{-\varphi}\geq\int_\Omega\langle
[F_\varphi+p(n-p)\lambda_\mathfrak{R}]g,g\rangle e^{-\varphi}\eeq
where $g\in C^\infty_p(\overline{\Omega})\cap{\rm
Dom}(d_\varphi^*).$ By H\"{o}rmander's density lemma(see [H1] or
[H2]), we know that \eqref{eq720} holds for any $g\in {\rm
Dom}(d_\varphi^*).$ Now the desired result follows from the estimate
\eqref{eq721} and lemma \ref{lem31}. For the regularity, we can
apply the procedure in section 5 to get the minimal solution in
$L^2_{p-1}(M,\varphi)$ and the interior regularity then follows from
the ellipticity of $d\delta_\varphi+\delta_\varphi d$.\eprf

\brem\label{rem71} By results in \cite{MM} and \cite{Se}, we know
that the curvature term
$$\sum_{i_1<\cdots<i_p}\langle\mathfrak{R}\xi^g_{i_1\cdots
i_p},\xi^g_{i_1\cdots i_p}\rangle\geq 0$$ for $p=2$ when $(M,ds^2)$
is assumed to have nonnegative complex sectional curvature(isotropic
sectional curvature when $n$ is even) . In this case, we have
instead of \eqref{eq721} the following apriori estimate
\beq\label{eq722}\int_\Omega(|dg|^2+|\delta_\varphi
g|^2)e^{-\varphi}\geq\int_\Omega\langle F_\varphi g,g\rangle
e^{-\varphi}\eeq for any $g\in C^\infty_2(\overline{\Omega})\cap{\rm
Dom}(d_\varphi^*).$ The same argument for Theorem \ref{thm71} also implies the following result:\\
Let $(M,ds^2)$ be a $n$-dimensional oriented $2$-convex Riemannian
manifold. Let $\varphi\in C^2(M)$ be a $2$-plurisubharmonic function
on $M$. If $(M,ds^2)$ has nonnegative complex sectional
curvature(isotropic sectional curvature when $n$ is even), then for
any closed $2$-form $f\in L^2_2(M,{\rm Loc})$ with
$$\int_M\langle F_\varphi^{-1}f,f\rangle e^{-\varphi}<\infty,$$there exists some
$1$-form $u \in L^2_1(M,\varphi)$ such that $$du=f  \ {\rm and} \
\int_M|u|^2e^{-\varphi}\leq \int_M\langle F_\varphi^{-1}f,f\rangle
e^{-\varphi}.$$Moreover, if $f$ and $\varphi$ are both assumed
additionally to be smooth then we can choose $u$ to be a smooth
form. \erem

As an easy corollary, we have the following result which is a
generalization of Theorem \ref{thm31} to Riemannian manifolds with
nonnegative curvature operator.

\bcor \label{cor71} Assume that $(M,ds^2)$ is $p$-convex and has
nonnegative curvature operator. Let $\varphi\in C^2(M)$ be a
$p$-plurisubharmonic function on $M$. Then for any closed $p$-form
$f\in L^2_p(M,{\rm Loc})$ with
$$\int_M\langle F_\varphi^{-1}f,f\rangle e^{-\varphi}<\infty,$$there exists some
$(p-1)$-form $u \in L^2_{p-1}(M,\varphi)$ such that $$du=f  \ {\rm
and} \ \int_M|u|^2e^{-\varphi}\leq \int_M\langle
F_\varphi^{-1}f,f\rangle e^{-\varphi}.$$Moreover, if $f$ and
$\varphi$ are both assumed to be smooth then we can choose $u$ to be
a smooth form. When $p=2$, it is enough to assume $(M,ds^2)$ has
nonnegative complex sectional curvature(isotropic sectional
curvature when $n$ is even).\ecor

\brem\label{rem72}All the results in sections 2-6 can be established
on Rimannian manifolds without any additional difficulty. For
Theorem \ref{thm21}, the minor difference is that the Levi-Civita
connection $D$ enters the derivatives and the gradient is taken
with respect to the underlying metric. To prove, on Riemannian manifolds,
these $L^2$-estimates obtained in sections 3-6, the only
modification is to use the estimate \eqref{eq721} to replace
\eqref{eq32}(or use Theorem \ref{thm71} to replace Theorem
\ref{thm31}). \erem

\sec{Geometric applications}

In this section, we will prove vanishing and finiteness theorems for
de Rham cohomology groups. The key is to control the curvature
term(in the basic estimate \eqref{eq721}) by choosing appropriate
weight functions. The main tool is a Carleman type estimate (Lemma
\ref{lem85}) which is uniform with respect to both of weights and domains. To
establish such an estimate, we will first prove a G{\aa}rding type
estimate(Lemma \ref{lem82}) which is also uniform w.r.t domains and
weights. Since the notion of $p$-convexity depends on the underlying
metric, we do not have the flexibility in the way of modifying the
metric as the complex analytic case(cf. \cite{AV} and \cite{D2}).\\

Solving $du=f$ in appropriate weighted $L^2$-space, we have the
following immediate corollary of Theorem \ref{thm71}.

\bpp \label{pp82}Let $(M,ds^2)$ be a  strictly $p$-convex
$n$-dimensional Riemannian manifold, $1\leq p\leq n$. Then for any
closed $f\in L^2_q(M,{\rm Loc})(p\leq q\leq n)$ there exists some
$(q-1)$-form $u\in L^2_{q-1}(M,{\rm Loc})$ such that $du=f$. In
particular, the de Rham cohomology group $H^{q}(M,\mathbb{R})=0$ for
every $p\leq q\leq n$.\epp

\bprfp{pp82}Since strict $p$-convexity implies
 strict $(p+1)$-convexity, it suffices to consider the case $q=p$. By using lemma \ref{lem81}
(i) with $\omega=p(n-p)\lambda_\mathfrak{R}$, one can find a
$p$-plurisubharmonic proper exhaustion function $\varphi\in
C^\infty(M)$ such that $F_\varphi+p(n-p)\lambda_\mathfrak{R}Id$ is
$p$-positive definite on $M$. Then
$\langle[F_\varphi+p(n-p)\lambda_\mathfrak{R}{\rm
Id}]^{-1}f,f\rangle$ is a continuous function on $M$. By Lemma
\ref{lem81} (iii), one can find a function $\psi\in C^\infty(M)$
such that $\psi-\varphi$ is $p$-plurisubharmonic and
$\int_M\langle[F_\varphi+p(n-p)\lambda_\mathfrak{R}{\rm
Id}]^{-1}f,f\rangle e^{-\psi}<\infty$. Consequently, we have
$$\int_M\langle[F_\psi+p(n-p)\lambda_\mathfrak{R}{\rm
Id}]^{-1}f,f\rangle
e^{-\psi}\leq\int_M\langle[F_\varphi+p(n-p)\lambda_\mathfrak{R}{\rm
Id}]^{-1}f,f\rangle e^{-\psi}<\infty.$$ It follows from Theorem
\ref{thm71} that there exists some $(p-1)$-form $u\in
L^2_{p-1}(M,\psi)$ such that $du=f$. To see the vanishing of
$H^{p}(M,\mathbb{R})$, it is enough to consider the minimal solution
of $du=f$ in $L^2_{p-1}(M,\psi)$ which is smooth provided $f\in
C^\infty_p(M)$. The proof is complete.\eprf

\brem\label{rem81}(i) As observed by Harvey and Lawson, proposition
\ref{pp82} also follows from Morse theory(see Theorem 4.16 in
\cite{HL2}).\\ (ii) By making an additional assumption on the
sectional curvature, we can prove the following vanishing result for
Riemannian manifolds which are strictly $p$-convex at infinity. Since $(M,ds^2)$
is strictly $p$-convex at infinity, $M$ can be exhausted by open subsets
with strictly $p$-convex boundary
$\Omega_1\Subset\Omega_2\Subset\cdots$. When $(M,ds^2)$ is assumed to have
nonnegative sectional curvature, by the main theorem in \cite{Sh},
we obtain $$H^q(\Omega_\nu,\mathbb{R})=0 \ {\rm for \ each } \ \nu\geq 1 \ {\rm and} \ 
p\leq q\leq n.$$ Taking the inverse limit implies that $$H^q(M,\mathbb{R})\cong \varprojlim
H^q(\Omega_\nu,\mathbb{R})=0 \ {\rm  for} \  p\leq q\leq n.$$ This is a generalization of Corollary \ref{cor51}. \erem

Combining the inequalities \eqref{eq92} and \eqref{eq93} below, we
will get a G{\aa}rding type estimate which is \textbf {uniform} with respect to both
$p$-convex domains $\Omega\Subset M$ and $p$-plurisubharmonic
weight functions $\varphi\in C^2(M)$ satisfying the condition \eqref{eq91}
below. The existence of such a weight is given by Lemma \ref{lem81} (ii).  In the sequel, we
will denote by $(d|_\Omega)^{*}_\varphi$ the adjoint of the maximal
differential operator $d|_\Omega: L^2_q(\Omega,\varphi)\rightarrow
L^2_{q+1}(\Omega,\varphi)$.

\blem \label{lem82} Let $M$ be an oriented $n$-dimensional manifold. Let $\varphi$ be a $C^2$ function which is
$p$-plurisubharmonic on $M$ and satisfies
\beq\label{eq91}F_\varphi+[\max_{p\leq \ell\leq
n}\ell(n-\ell)\lambda_\mathfrak{R}-1]{\rm Id} \ {\rm is} \ p{\rm
-positive \ outside\ a \ compact \ subset} \ S\subseteq M.\eeq For
any bounded open set $\Omega$ with $p$-convex boundary and any open
neighborhoods $U\Subset U_1\subseteq \Omega$ of $S$ in $\Omega$,
there is a constant $A=A(S,U,U_1)>0$ such that
\beq\label{eq92}\int_\Omega \Big(|dg|^2+|(d|_\Omega)_\varphi^{*}
g|^2+|g|^2\Big)e^{-\varphi}\geq A\int_U |Dg|^2e^{-\varphi}\eeq and
\beq\label{eq93}\int_\Omega\Big(|dg|^2+|(d|_\Omega)_\varphi^{*}
g|^2\Big)e^{-\varphi}+\int_U|g|^2e^{-\varphi}\geq A\int_\Omega
|g|^2e^{-\varphi}\eeq hold for every $g\in {\rm
Dom}(d|_{\Omega})_\varphi^{*}\cap{\rm Dom}(d|_{\Omega})\subseteq
L^2_q(\Omega,\varphi),p\leq q\leq n.$ \elem

\bprfl{lem82}Let $g\in C_q^\infty(\overline{\Omega}), p\leq q\leq
n$. Choosing a cut-off function $\chi_1\in C^\infty(\Omega)$ such
that
$$\chi_1|_{U}\equiv1 \ {\rm and} \ {\rm Supp}\chi_1\subseteq U_1.$$It follows
from \eqref{eq716} and \eqref{eq719} that\beq\int_{\Omega}
\Big(|d(\chi_1g)|^2+|\delta_\varphi
(\chi_1g)|^2\Big)e^{-\varphi}&\geq&\int_\Omega
\Big(|D(\chi_1g)|^2-q(n-q)A_1|(\chi_1g)|^2\Big)e^{-\varphi}\nno\\&\geq&\int_{U}
|Dg|^2e^{-\varphi}-q(n-q)A_1\int_{\Omega}
|\chi_1g|^2e^{-\varphi}\nno\eeq where $A_1>0$ is a constant such
that $\lambda_\mathfrak{R}\geq -A_1$ on $U_1$. Therefore, we obtain
\beq\label{eq8103}\int_{\Omega} \Big\{|dg|^2+|\delta_\varphi
g|^2+A_2^2[q(n-q)A_1+2]|g|^2\Big\}e^{-\varphi}\geq\frac{1}{2}\int_{U}
|Dg|^2e^{-\varphi}\eeq where
$A_2:=\sup_\Omega(\frac{|\chi_1|}{\sqrt{2}}+|\nabla\chi_1|)$ and
$g\in
C^\infty_q(\overline{\Omega})$.\\
Let $\chi_2$ be a smooth function on $\Omega$ satisfying
$$\chi_2|_S\equiv 0 \ {\rm and} \ \chi_2|_{\Omega\setminus U}\equiv
1.$$ Set $A_3:= \sup_\Omega|\nabla\chi_2|$. For any $g\in
C^\infty_q(\overline{\Omega})\cap {\rm
Dom}((d|_{\Omega})_\varphi^{*})$, by using \eqref{eq716} and
\eqref{eq719} again, we have\beq\int_{\Omega}\Big(|d(\chi_2
g)|^2+|\delta_\varphi (\chi_2
g)|^2\Big)e^{-\varphi}&\geq&\int_\Omega\Big\langle(F_\varphi+q(n-q)\lambda_\mathfrak{R}{\rm
Id})\chi_2
g,\chi_2 g\Big\rangle e^{-\varphi}\nno\\
&\geq&\int_\Omega |\chi_2 g|^2e^{-\varphi}\nno\\
&\geq&\int_{\Omega\setminus U} |g|^2e^{-\varphi}\nno\eeq which
implies that\beq\int_{\Omega} \Big(|dg|^2+|\delta_\varphi
g|^2\Big)e^{-\varphi}+2A_3^2\int_{U}|g|^2e^{-\varphi}\geq\frac{1}{2}\int_{\Omega\setminus
U}|g|^2e^{-\varphi}\nno\eeq and consequently,\beq\label{eq8104}
\int_{\Omega} \Big(|dg|^2+|\delta_\varphi
g|^2\Big)e^{-\varphi}+(2A_3^2+\frac{1}{2})\int_{U}|g|^2e^{-\varphi}\geq\frac{1}{2}\int_\Omega
|g|^2e^{-\varphi}.\eeq By H\"{o}mander's density
lemma(cf.\cite{H1} or \cite{H2}), the estimates \eqref{eq8103} and
\eqref{eq8104} are both valid for $g\in {\rm
Dom}(d|_{\Omega})_\varphi^{*}\cap{\rm Dom}(d|_{\Omega})\subseteq
L^2_q(\Omega,\varphi),p\leq q\leq n.$ \eprf

By a compactness argument, the next result follows from Lemma
\ref{lem82}.

\blem \label{lem83} Let $\varphi\in C^2(M)$ be a
$p$-plurisubharmonic function satisfying \eqref{eq91}. For any
bounded open set $\Omega$ with $p$-convex boundary which contains
the subset $S$ in \eqref{eq91}, ${\rm
Ker}(d|_\Omega)_\varphi^{*}\cap{\rm Ker}(d|_\Omega)\subseteq
L^2_q(\Omega,\varphi)$ is finite dimensional  and we have the
orthogonal decomposition \beq\label{eq8106}{\rm
Ker}(d|_\Omega)=\Big({\rm Ker}(d|_\Omega)_\varphi^{*}\cap{\rm
Ker}(d|_\Omega)\Big)\oplus{\rm Im}(d|_\Omega)\subseteq
L^2_q(\Omega,\varphi),p\leq q\leq n.\eeq\elem

\bprfl{lem83}Fix open neighborhoods $U\Subset U_1\subseteq \Omega$
of $S$ in $\Omega$ such that $U$ has smooth boundary. Let
$\{g_\nu\}\subseteq {\rm Ker}(d|_\Omega)_\varphi^{*}\cap{\rm
Ker}(d|_\Omega)$ be a sequence of $q$-forms with $\|g_\nu\|_\varphi$
bounded and $\|dg_\nu\|_\varphi\rightarrow 0,
\|(d|_\Omega)_\varphi^{*} g_\nu\|_\varphi\rightarrow 0.$ In view of
\eqref{eq92} and the Rellich-Kondrakov theorem, we can pass to a
subsequence and thereby assume that $\{g_\nu|_{\Omega_1}\}$
converges in $L^2_q(U,\varphi)$. On the other hand, \eqref{eq93}
implies that $\{g_\nu$ is a cauchy sequence in
$L^2_q(\Omega,\varphi)$. Therefore, there exists a $g\in
L^2_q(\Omega,\varphi)$ such that $g_\nu\rightarrow g$ in
$L^2_q(\Omega,\varphi)$.By applying Lemma \ref{lem84} below to the
weighted $L^2$-de Rham complex $$\cdots\rightarrow
L_{q-1}^2(\Omega,\varphi)\overset{T=d|_\Omega}{\rightarrow}
L_{q}^2(\Omega,\varphi)\overset{S=d|_\Omega}{\rightarrow}L_{q}^2(\Omega,\varphi)\rightarrow\cdots$$
we get the desired results.\eprf

In the proof of Lemma \ref{lem83}, we have used the following
result.

\blem \label{lem84}(Theorems 1.12 and 1.13 in \cite{H1}) Let
$H_1\sr{T}{\to}H_2\sr{S}{\to}H_3$ be a complex of closed and densely
defined operators between Hilbert spaces. Assume that from every
sequence $g_\nu\in {\rm Dom}(T^{*})\cap{\rm Dom}(S)$ with
$\|g_\nu\|_{H_2}$ bounded and $T^{*}g_\nu\rightarrow 0$ in $H_1$,
$Sg_\nu\rightarrow 0$ in $H_3$, one can select a strongly convergent
subsequence. Then there exists a constant $C>0$ such that
\beq\label{eq8100}\|g\|_{H_2}^2\le
C^2(\|T^{*}g\|_{H_1}^2+\|Sg\|_{H_3}^2)\eeq holds for any $g\in {\rm
Dom}(T^{*})\cap {\rm Dom}(S)\cap ({\rm Ker}T^{*}\cap{\rm
Ker}S)^\perp$ and ${\rm Ker}T^{*}\cap{\rm Ker}S$ is finite
dimensional. Moreover, when the above estimate \eqref{eq8100} holds,
we also the following orthogonal decomposition
\beq\label{eq8101}{\rm Ker}S=({\rm Ker}T^{*}\cap{\rm
Ker}S)\oplus{\rm Im}T.\eeq \elem

\brem\label{rem82}Since $L^2_{*}(\Omega,\varphi) =L^2_{*}(\Omega)$,
in the orthogonal decomposition \eqref{eq8106}, the left hand side
and the second summand on the right hand side are independent of the
choice of $\varphi$. Different choices of $\varphi$ result in
different complementary subspaces of ${\rm Im}(d|_\Omega)$ in ${\rm
Ker}(d|_\Omega)$. \erem

We can deduce from Lemma \ref{lem82} a Carleman type inequality
which is uniform with respect to an increasing sequence of open subsets and
weight functions. To formulate such estimates, we introduce an
increasing sequence of convex increasing functions $\chi_{\nu}\in
C^2(\mathbb{R}),\nu=1,2,\cdots$ such that
\beq\label{eq1}\chi_{\nu}(t)\equiv 0 \ {\rm for} \ t\leq 0 \ {\rm
and} \ \nu=1,2,\cdots, \lim_{\nu\rightarrow
+\infty}\chi_{\nu}(t)=+\infty \ {\rm for} \ t>0.\eeq

\blem \label{lem85} Let $\varphi\in C^2(M)$ be a
$p$-plurisubharmonic function satisfying \eqref{eq91}. Assume that
the subset $S$ in \eqref{eq91} is contained in $U:=\{x\in M \ | \
\varphi(x)<0\}$ and that $U$ has smooth boundary. Then for any
sequence $\Omega_1\Subset\Omega_2\Subset\cdots\subseteq M$ of smooth
open subsets with $p$-convex boundary such that
\beq\label{eq1234}U\Subset\displaystyle{\cup_{\nu\geq
1}}\Omega_\nu,\eeq there exist constants $C>0$ and $\nu_0>0$ such
that
\beq\label{eq2}\int_{\Omega_\mu}|f|^2e^{-\varphi-\chi_\nu\circ\varphi}\leq
C^2\int_{\Omega_\mu}\Big(|(d|_{\Omega_\mu})_{\varphi+\chi_\nu\circ\varphi}^{*}f|^2
+|df|^2\Big)e^{-\varphi-\chi_\nu\circ\varphi}\eeq for every
$\mu,\nu\geq\nu_0$ and every $f\in{\rm
Dom}((d|_{\Omega_\mu})_{\varphi+\chi_\nu\circ\varphi}^{*})\cap{\rm
Dom}(d|_{\Omega_\mu})\subseteq
L^2_q(\Omega_\mu,{\varphi+\chi_\nu\circ\varphi})$ which satisfies
\beq\label{eq3}\int_U\langle f,g\rangle e^{-\varphi}=0, \forall
g\in{\rm Ker}((d|_U)_\varphi^{*})\cap{\rm Ker}(d|_U)\subseteq
L^2_q(U,\varphi)\eeq where $p\leq q\leq n$ and $\{\chi_\nu\}$ is any
increasing sequence consists of convex increasing functions
satisfying \eqref{eq1}. \elem

\bprfl{lem85}We proceed by contradiction. Since
$U\Subset\displaystyle{\cup_{\nu\geq 1}}\Omega_\nu$, we can assume,
without loss of generality, that $U\Subset\Omega_1.$

It is easy to see that each $\varphi+\chi_\nu\circ\varphi(\nu\geq
1)$ satisfies the condition \eqref{eq91} with the same subset $S$.
By Lemma \ref{lem82}, we know that \eqref{eq92} and \eqref{eq93}
hold for all subsets $\Omega_\mu$ and weight functions
$\varphi+\chi_\nu\circ\varphi(\mu,\nu\geq 1)$. It is easy to see, by
fixing an open set $U_1$ such that $U\Subset U_1\subseteq \Omega_1$,
that the constant $A$ in Lemma \ref{lem82} is independent of
$\mu,\nu\geq 1$.

If the conclusion were false, by passing to subsequences of
$\{\Omega_\mu\}_{\mu\geq 1}$ and $\{\chi_\nu\}_{\nu\geq 1}$(as the
conditions \eqref{eq1} and \eqref{eq1234} are both fulfilled for any
subsequence), we may assume that there exists a sequence of
$f_\nu\in{\rm
Dom}((d|_{\Omega_\nu})_{\varphi+\chi_\nu\circ\varphi}^{*})\cap{\rm
Dom}(d|_{\Omega_\nu})\subseteq
L^2_q(\Omega_\nu,{\varphi+\chi_\nu\circ\varphi})(\nu\geq 1)$ with
the following
properties\beq\label{eq4}\int_{\Omega_\nu}|f_\nu|^2e^{-\varphi-\chi_\nu\circ\varphi}=1,\eeq
\beq\label{eq5}\int_{\Omega_\nu}\Big(|(d|_{\Omega_\nu})_{\varphi+\chi_\nu\circ\varphi}^{*}f_\nu|^2
+|df_\nu|^2\Big)e^{-\varphi-\chi_\nu\circ\varphi}\leq \nu^{-1},\eeq
\beq\label{eq6} \int_U\langle f_\nu,g\rangle e^{-\varphi}=0, \forall
g\in{\rm Ker}((d|_U)_\varphi^{*})\cap{\rm Ker}(d|_U)\subseteq
L^2_q(U,\varphi).\eeq

By \eqref{eq92}, \eqref{eq4} and \eqref{eq5}, we get
\beq\int_U|Df_\nu|^2e^{-\varphi}&=&
\int_U|Df_\nu|^2e^{-\varphi-\chi_\nu\circ\varphi}\nno \\
&\leq&
A^{-1}\int_{\Omega_\nu}\Big(|f_\nu|^2+|(d|_{\Omega_\nu})_{\varphi+\chi_\nu\circ\varphi}^{*}f_\nu|^2
+|df_\nu|^2\Big)e^{-\varphi-\chi_\nu\circ\varphi}\nno \\ &\leq&
A^{-1}(1+\nu^{-1})\nno\eeq and \beq\int_U|f_\nu|^2e^{-\varphi}=
\int_U|f_\nu|^2e^{-\varphi-\chi_\nu\circ\varphi}\leq
\int_{\Omega_\nu}|f_\nu|^2e^{-\varphi-\chi_\nu\circ\varphi}=1.\nno\eeq
The Rellich-Kondrakov theorem implies that we may assume, by passing
to a subsequence, that
\beq\label{eq7}\lim_{\nu\rightarrow+\infty}f_\nu=f \ {\rm in} \
L^2_q(U,\varphi).\eeq Taking into account of \eqref{eq5}, we also have
\beq\lim_{\nu\rightarrow+\infty}df_\nu=0 \ {\rm in} \
L^2_q(U,\varphi)\nno\eeq which implies that
\beq\label{eq101}f\in{\rm Ker}(d|_U)\subseteq L^2_q(U,\varphi).\eeq
Taking limit in \eqref{eq6}, we obtain\beq\label{eq8}\int_U\langle
f,g\rangle e^{-\varphi}=0, \forall g\in{\rm
Ker}((d|_U)_\varphi^{*})\cap{\rm Ker}(d|_U)\subseteq
L^2_q(U,\varphi).\eeq  

\par

Now set $$g_\nu=e^{-\chi_\nu\circ\varphi}f_\nu(\nu\geq 1),$$ by using
\eqref{eq4} and \eqref{eq5} respectively, we
have\beq\label{eq9}\int_{\Omega_1}|g_\nu|^2e^{-\varphi+\chi_\mu\circ\varphi}\leq\int_{\Omega_\nu}|g_\nu|^2e^{-\varphi+\chi_\nu\circ\varphi}
=\int_{\Omega_\nu}|f_\nu|^2e^{-\varphi-\chi_\nu\circ\varphi}=1,\eeq
and
\beq\label{eq10}\int_{\Omega_1}|(d|_{\Omega_\nu})_{\varphi}^{*}g_\nu|^2
e^{-\varphi}&\leq&\int_{\Omega_\nu}|(d|_{\Omega_\nu})_{\varphi}^{*}g_\nu|^2
e^{-\varphi+\chi_\nu\circ\varphi}\nno \\&=&
\int_{\Omega_\nu}|(d|_{\Omega_\nu})_{\varphi+\chi_\nu\circ\varphi}^{*}f_\nu|^2
e^{-\varphi-\chi_\nu\circ\varphi}\leq\nu^{-1}\eeq for any
$\nu\geq\mu\geq 1$. By \eqref{eq9}, we may assume
\beq\label{eq102}g_\nu\overset{{\rm w}}{\rightharpoonup} g \ {\rm
in}\ L^2_q(\Omega_1,\varphi-\chi_\mu\circ\varphi)\ {\rm as} \
\nu\rightarrow+\infty\eeq for any $\mu\geq 1$. Combining \eqref{eq9}
and \eqref{eq102} gives
\beq\int_{\Omega_1}|g|^2e^{-\varphi+\chi_\mu\circ\varphi}\leq
1\nno\eeq for any $\mu\geq 1$, which implies that
\beq\label{eq11113}{\rm Supp}g\subseteq U.\eeq By \eqref{eq102}, we
know that \beq\label{eq1110}g_\nu|_{\Omega_1}\rightarrow g \ {\rm in
\ the \ sense \ of \ distribution}\eeq as $\nu\rightarrow+\infty$.
Consequently, \beq\label{eq104}\delta_\varphi
g_\nu|_{\Omega_1}\rightarrow \delta_\varphi g\eeq in the sense of
distribution, as $\nu\rightarrow+\infty$. Meanwhile, we know by
\eqref{eq10} that \beq\label{eq105}\delta_\varphi
g_\nu|_{\Omega_1}\rightarrow 0 \ {\rm in}\ L^2_q(\Omega_1,\varphi)\
{\rm as} \ \nu\rightarrow+\infty.\eeq Combining \eqref{eq104} and
\eqref{eq105}, we get \beq\label{eq106}\delta_\varphi g=0 \ {\rm on}
\ \Omega_1\eeq in the sense of distribution. From \eqref{eq11113}
and \eqref{eq106}, it follows that \beq\label{eq107}g|_U\in {\rm
Ker}((d|_U)^{*}_\varphi).\eeq

By the definition of $g_\nu$, we know $$g_\nu=f_\nu \ {\rm on} \
U$$ which, together with  \eqref{eq7}, \eqref{eq101} and
\eqref{eq1110}, implies that \beq\label{eq108}g|_U=f\in{\rm
Ker}(d|_U)\subseteq L^2_q(U,\varphi).\eeq From
\eqref{eq8},\eqref{eq107} and \eqref{eq108}, it follows that
\beq\label{eq109}\lim_{\nu\rightarrow+\infty}f_\nu|_U=f=0 \ {\rm in}
\ L^2_q(U,\varphi).\eeq On the other hand, by \eqref{eq93},
\eqref{eq4} and \eqref{eq5}, we have
\beq\nu^{-1}+\int_U|f_\nu|^2e^{-\varphi}
&\geq&\int_{\Omega_\nu}\Big(|(d|_{\Omega_\nu})_{\varphi+\chi_\nu\circ\varphi}^{*}f_\nu|^2
+|df_\nu|^2\Big)e^{-\varphi-\chi_\nu\circ\varphi}+\int_U|f_\nu|^2e^{-\varphi}\nno
\\ &\geq&
A\int_{\Omega_\nu}|f_\nu|^2e^{-\varphi-\chi_\nu\circ\varphi}=A\nno\eeq
Letting $\nu\rightarrow+\infty$ and using \eqref{eq109}, we get the
contradiction $0\geq A$ which completes the proof.\eprf

\bthm\label{thm81} Let $(M,ds^2)$ be an $n$-dimensional
Riemannian manifold which is strictly $p$-convex at infinity($1\leq
p\leq n$). Then the de Rham cohomology group $H^{q}(M,\mathbb{R})$
is finite dimensional for every $p\leq q\leq n$.\ethm

\bprft{thm81}Let $\pi : \widetilde{M}\to M$ be the orientation covering of $M$, then we know by definition that $\widetilde{M}$,  endowed with the pulled back metric $\pi^*ds^2$, is again strictly $p$-convex at infinity. Since $\pi : \widetilde{M}\to M$ is a double covering, the induced homomorphism $\pi^*:H^{q}(M,\mathbb{R})\to H^{q}(\widetilde{M},\mathbb{R})$ is injective for every $q$. By passing to $\widetilde{M}$, we may assume without loss of generality that $M$ is oriented. We will deduce Theorem \ref{thm81} as a consequence of Lemmas \ref{lem83} and \ref{lem85} in the following three steps.  

\par

Step 1. By Lemma \ref{lem81} (ii), there is a proper
exhaustion function $\varphi\in C^\infty(M)$ satisfying the
hypotheses of Lemma \ref{lem85} where $\Omega_\nu:=\{x\in M \ | \
\varphi(x)<\nu\}, \nu=1,2,\cdots.$ From Lemma \ref{lem83}(choose
$\Omega$ to be the subset $U$ in Lemma \ref{lem85}), it is
sufficient to prove that the natural homomorphism from
the de Rham cohomology $H^q(M,\mathbb{R})$ to the $L^2$-cohomology $L^2H^{q}(U):=\frac{\{f\in L^2_q(U) \ | \
df=0\}}{\{du\in L^2_q(U) \ | \ u\in L^2_{q-1}(U\}}$, given by the
pullback of the inclusion map, is injective for any $p\leq q\leq n$.

Step 2. By Corollary \ref{cor31}, we have the following fine
resolution of the constant sheaf $\mathbb{R}$
$$0\rightarrow\mathbb{R}\rightarrow
\mathscr{A}_0\overset{d}{\rightarrow}
\mathscr{A}_1\overset{d}{\rightarrow}
\mathscr{A}_2\overset{d}{\rightarrow}\cdots$$where
$\mathscr{A}_q(V):=\{f\in L^2_q(V,{\rm Loc}) \ | \ df\in
L^2_{q+1}(V,{\rm Loc})\}$ for any open subset $V\subseteq M$ and
$0\leq q \leq n$. Hence it suffices to show that for any closed
$q$-form $h\in L^2_q(M,{\rm Loc})$ if \beq\label{eq1111}h|_U=du \
{\rm where} \ u\in L^2_{q-1}(U)\eeq then there exists a $(q-1$)-form
$\widetilde{u}\in L^2_{q-1}(M,{\rm Loc})$ such that
$d\widetilde{u}=h$ holds on $M$ in the sense of distribution.

Step 3. By Lemma \ref{lem81} (iii), one can find some function
$\psi\in C^\infty(M)$ such that $\varphi\equiv\psi$ when
$\varphi<1$, $\psi-\varphi$ is $p$-plurisubharmonic and that
\beq\label{eq1114}\int_M|h|^2e^{-\psi}<+\infty.\eeq It is easy to
see that $\psi$ still satisfies the hypotheses of Lemma \ref{lem85}
with the same $S$ and $U$. By Lemma \ref{lem85}, there are constants
$C>0$ and $\nu_0>0$ such that
\beq\int_{\Omega_\nu}|f|^2e^{-\psi-\chi_{\nu_0}\circ\psi}\leq
C^2\int_{\Omega_\nu}\Big(|(d|_{\Omega_\nu})_{\psi+\chi_{\nu_0}\circ\psi}^{*}f|^2
+|df|^2\Big)e^{-\psi-\chi_{\nu_0}\circ\psi}\nno\eeq holds for every
$f\in{\rm
Dom}((d|_{\Omega_\nu})_{\psi+\chi_{\nu_0}\circ\psi}^{*})\cap{\rm
Dom}(d|_{\Omega_\nu})\subseteq
L^2_q(\Omega_\nu,{\psi+\chi_{\nu_0}\circ\psi})$ satisfying
\eqref{eq3} where $\nu=\nu_0,\nu_0+1,\nu_0+2,\cdots$. 

From the above estimate and Lemma\ref{lem31}, we know that for any closed $q$-form $h_\nu\in
L^2_q(\Omega_\nu,{\psi+\chi_{\nu_0}\circ\psi})$ satisfying
\eqref{eq3}, there exists, for each $\nu\geq\nu_0$, some
$(q-1)$-form $u_\nu\in
L^2_{q-1}(\Omega_\nu,\psi+\chi_{\nu_0}\circ\psi)$ such that
\beq\label{eq1113}du_\nu=h_\nu,
\int_{\Omega_\nu}|u_\nu|^2e^{-\psi-\chi_{\nu_0}\circ\psi}\leq
C^2\int_{\Omega_\nu}|h_\nu|^2e^{-\psi-\chi_{\nu_0}\circ\psi}\leq
C^2\int_{\Omega_\nu}|h_\nu|^2e^{-\psi}\eeq By
\eqref{eq1111},\eqref{eq1114} and \eqref{eq1113}, we get some
$u_\nu\in L^2_{q-1}(\Omega_\nu,\psi+\chi_{\nu_0}\circ\psi)$ such
that \beq\label{eq1112}du_\nu=h|_{\Omega_\nu},
\int_{\Omega_\nu}|u_\nu|^2e^{-\psi-\chi_{\nu_0}\circ\psi}\leq
C^2\int_{\Omega_\nu}|h|^2e^{-\psi}<+\infty\eeq for each
$\nu\geq\nu_0$. Now the desired solution $\widetilde{u}\in
L^2_{q-1}(M,\psi+\chi_{\nu_0}\circ\psi)\subseteq L^2_{q-1}(M,{\rm
Loc})$ follows from using \eqref{eq1112} to take weak limits.\eprf

As an intermediate consequence, we have

\bcor\label{cor82} Let $(M,ds^2)$ be a oriented $n$-dimensional
Riemannian manifold, and $\Omega\Subset M$ be an open subset with
strictly $p$-convex boundary, then the de Rham cohomology group
$H^q(\Omega,\mathbb{R})$ is finite dimensional for every $p\leq
q\leq n.$\ecor

\bprfc{cor82}By Lemma 3.17 in \cite{HL3}, we know that $\Omega$ is
strictly $p$-convex at infinity w.r.t the induced metric from $M$.
Thus the finiteness result follows from Theorem \ref{thm81}.\eprf

In the above proof of Theorem \ref{thm81}, Lemma \ref{lem85} is
applied to a fixed weight function and a sequence of domains. If we
apply Lemma \ref{lem85} to a fixed domain and a sequence of weight
functions, then we achieve the following approximation result.

\bthm\label{thm82}Let $\varphi\in C^2(M)$ be a $p$-plurisubharmonic
function satisfying \eqref{eq91}. Assume that the subset $S$ in
\eqref{eq91} is contained in $U:=\{x\in M \ | \ \varphi(x)<0\}$ and
that $U
$ has smooth boundary. Let $\Omega\Subset M$ be an open
subset with $p$-convex boundary such that $U\Subset\Omega$. Then for
any closed $(q-1)$-form $u\in L^2_{q-1}(U)$ there exists a sequence
of closed $(q-1)$-forms $u_\nu\in L^2_{q-1}(\Omega)$ such that
$u_{\nu}|_U\rightarrow u \ {\rm in} \ L^2_{q-1}(U)$ where $p\leq
q\leq n$. \ethm

\bprft{thm82}It is easy to see that $${\rm Ker}(d|_\Omega)\subseteq
{\rm Ker}(d|_U)\subseteq L^2_{q-1}(U).$$ The desired conclusion is
$\overline{{\rm Ker}(d|_\Omega)}\supseteq {\rm Ker}(d|_U)$ where
$\overline{\cdot}$ means taking the closure in $L^2_{q-1}(U)$. Since
${\rm Ker}(d|_U)\subseteq L^2_{q-1}(U)$ is closed, it suffices to
show
$${\rm Ker}(d|_\Omega)^{\perp}\subseteq {\rm Ker}(d|_U)^{\perp}$$ where $\cdot^{\perp}$ means taking
the orthogonal complement in the Hilbert space $L^2_{q-1}(U)$.

Let $u\in {\rm Ker}(d|_\Omega)^{\perp}\subseteq L^2_{q-1}(U)$, we
extend $u$ to be an element $\widetilde{u}\in L^2_{q-1}(\Omega)$ by
setting $\widetilde{u}=0$ outside $U$. The condition $u\in {\rm
Ker}(d|_\Omega)^{\perp}\subseteq L^2_{q-1}(U)$ implies that
$\widetilde{u}$ lies in the orthogonal complement of ${\rm
Ker}(d|_\Omega)$ in $L^2_{q-1}(\Omega)$. Let $\{\chi_\nu\}$ be an
increasing sequence consists of convex increasing functions
satisfying \eqref{eq1}, then we know that
\beq\widetilde{u}e^{\varphi+\chi_\nu\circ\varphi}\in {\rm
Ker}(d|_\Omega)^{\perp} \subseteq
L^2_{q-1}(\Omega,\varphi+\chi_\nu\circ\varphi),\nno\eeq as before, 
$\cdot^{\perp}$ means taking the orthogonal complement in the Hilbert space
$L^2_{q-1}(\Omega,\varphi+\chi_\nu\circ\varphi)$. By Lemma
\ref{lem83}, it follows that $${\rm Ker}(d|_\Omega)^{\perp}={\rm
Im}((d|_\Omega)^{*}_{\varphi+\chi_\nu\circ\varphi})\subseteq
L^2_{q-1}(\Omega,\varphi+\chi_\nu\circ\varphi).$$ Hence we can find
a unique $f_\nu\in {\rm
Dom}((d|_\Omega)^{*}_{\varphi+\chi_\nu\circ\varphi})\cap {\rm
Ker}((d|_\Omega)^{*}_{\varphi+\chi_\nu\circ\varphi})^{\perp}\subseteq
L^2_q(\Omega,\varphi+\chi_\nu\circ\varphi)$
 such that \beq\label{eq2222}(d|_\Omega)^{*}_{\varphi+\chi_\nu\circ\varphi}f_\nu=
e^{\varphi+\chi_\nu\circ\varphi}\widetilde{u} \ {\rm on} \
\Omega\eeq for each $\nu\geq 1$.

Since ${\rm
Ker}((d|_\Omega)^{*}_{\varphi+\chi_\nu\circ\varphi})^{\perp}({\rm
in} \ L^2_{q-1}(\Omega,\varphi+\chi_\nu\circ\varphi))\subseteq {\rm
Ker}((d|_U)^{*}_\varphi)^{\perp}({\rm in} \ L^2_{q-1}(U,\varphi))$,
by Lemma \ref{lem85} and \eqref{eq2222}, there are constants $C>0$
and $\nu_0>0$ such
that\beq\label{eq3333}\int_\Omega|f_\nu|^2e^{-\varphi-\chi_\nu\circ\varphi}&\leq&
C^2\int_\Omega|e^{\varphi+\chi_\nu\circ\varphi}\widetilde{u}|^2e^{-\varphi-\chi_\nu\circ\varphi}\nno
\\ &=&
C^2\int_U|u|^2e^\varphi\eeq for any $\nu\geq\nu_0.$  Set
$g_\nu=e^{-\varphi-\chi_\nu\circ\varphi}f_\nu(\nu\geq 1)$, then by
\eqref{eq2222} and \eqref{eq3333} we
get\beq\label{eq4444}(d|_\Omega)^{*}g_\nu=\widetilde{u}\eeq and
\beq\label{eq5555}\int_\Omega|g_\nu|^2e^{\varphi+\chi_\nu\circ\varphi}\leq
C^2\int_U|u|^2e^\varphi<+\infty\eeq for any $\nu\geq\nu_0.$ Estimate
\eqref{eq5555} implies that
\beq\label{eq7777}\int_\Omega|g_\nu|^2e^{\varphi+\chi_\mu\circ\varphi}\leq
C^2\int_U|u|^2e^\varphi<+\infty\eeq for any $\nu\geq\mu\geq\nu_0$.
Therefore there is a weak limit, denoted by $g$,  of $g_\nu$ in
$L^2_{q-1}(\Omega,\varphi+\chi_\mu\circ\varphi)$ for any
$\mu\geq\nu_0$(note that $g$ is independent of $\mu$). It follows
from \eqref{eq4444} that \beq\label{eq8888}\delta g=\widetilde{u} \
{\rm on} \ \Omega\eeq in the sense of distribution. Letting
$\nu\rightarrow+\infty$, \eqref{eq7777}
gives\beq\int_\Omega|g|^2e^{\varphi+\chi_\mu\circ\varphi}\leq
C^2\int_U|u|^2e^\varphi<+\infty\nno\eeq for any $\mu\geq\nu_0$.
Taking limit as $\mu\rightarrow+\infty$ yields\beq\label{eq6666}{\rm
Supp g}\subseteq U.\eeq Combining \eqref{eq8888} and \eqref{eq6666}
shows that $g\in {\rm Dom}((d|_U)^{*})$ and consequently $u\in {\rm
Im}((d|_U)^{*})={\rm Ker}(d|_U)^{\perp}\subseteq L^2_{q-1}(U)$. The
proof is complete.\eprf

\addcontentsline{toc}{section}{Acknowledgements}
{\bf{Acknowledgements}}: We would like to thank Professor Paul Yang
for valuable suggestions which enabled us to improve an earlier
version of Theorem \ref{thm81}.  The authors are grateful to referees for very careful reading helpful suggestions.

\addcontentsline{toc}{section}{References}

\end{document}